\documentclass[12pt]{amsart}
\usepackage{amssymb}
\usepackage[active]{srcltx}
\usepackage{epsfig}
\usepackage{amscd}
\usepackage{a4wide}
\usepackage{amsthm}
\usepackage{tikz}
\usepackage{tikz-cd} 
\usetikzlibrary{matrix}

\newtheorem{theorem}{Theorem}[section]
\newtheorem{lemma}[theorem]{Lemma}
\newtheorem{proposition}[theorem]{Proposition}
\newtheorem{corollary}[theorem]{Corollary}

\theoremstyle{definition}
\newtheorem{definition}[theorem]{Definition}
\newtheorem{example}[theorem]{Example}
\newtheorem{remark}[theorem]{Remark}

\numberwithin{equation}{section}

\DeclareMathOperator{\Sphere}{S}
\newcommand{\bc}{{\mathbb C}}

\begin{document}
\title [Symplectic quandle Method and $SL(2,\bc)$-representations of 2-bridge knots] {Symplectic quandle Method and  $SL(2,\bc)$-representations of 2-bridge Knots}
\author{Kyeonghee Jo and Hyuk Kim}
% \footnotetext{}
  \address{Division of Liberal Arts and sciences,
  Mokpo National Maritime University, Mokpo, Chonnam, 530-729, Korea  }
\email{khjo@mmu.ac.kr}
\address{Department of Mathematical Sciences, Seoul National University, Seoul, 08826, Korea}
\email{hyukkim@snu.ac.kr}
%\thanks {}
\maketitle
\begin{abstract}
In this paper, we extend the symplectic quandle method, previously employed in our study of parabolic representations of knot groups, to investigate the general $SL(2,\mathbb{C})$-representations of 2-bridge ``kmot" groups.
%In this paper, we study $SL(2,\bc)$-representations of 2-bridge ``kmot" groups by extending the symplectic quandle method that we used in studying parabolic representations in the previous paper. 
We introduce a `generalized symplectic quandle structure' corresponding to  ($\mathcal{D}_M$, conjugation) for each $M\in\mathbb C\setminus \{0,1,-1\}$,  where $\mathcal{D}_M=\{A\in SL(2,\bc)\mid tr(A)= M+M^{-1} \}$.  By converting the system of conjugation quandle equations to that of generalized symplectic quandle equations, we obtain a simpler expression for the 2-variable Riley polynomial and derive some recursive formulas for Riley polynomials and Alexander polynomials. This approach enables us to effectively compute the A-polynomials, allowing us to obtain numerous previously unknown A-polynomials within minutes using Mathematica. 
\end{abstract}

\tableofcontents

\section{Introduction}
A 2-bridge ``kmot" group, denoted as $G$, is a group with the following presentation:
  $$G=\langle x,y \,|\,wx=yw \rangle,$$ where $w$ is of the form
\begin{equation}
w=x^{\epsilon_1}y^{\epsilon_2}x^{\epsilon_3}y^{\epsilon_4}\cdots x^{\epsilon_{\alpha-2}}y^{\epsilon_{\alpha-1}} \,\,(\alpha \,\text{is an odd integer}\,\geq 3)
\end{equation}
and  $\epsilon_i=\epsilon_{\alpha-i}=\pm 1$ for $i=1, \cdots, \alpha-1$. It is well known that all the non-abelian $SL(2,\bc)$-representations of 2-bridge kmot groups are given by the solutions of the 2-variable Riley polynomial  
$$\mathcal{R}(\lambda,s)=0, s=M^2+\frac{1}{M^2}-2$$ 
in $\mathbb Z[\lambda,s]$ (see \cite{Riley2}), and each solution corresponds to a $SL(2,\bc)$-representation $\rho$ such  that
\begin{equation*}
\rho(x)=\begin{pmatrix}
M & 1\\
  0 & \frac{1}{M} 
 \end{pmatrix} \quad \text{and}\quad
\rho(y)=\begin{pmatrix}
M & 0\\
  \lambda& \frac{1}{M} 
 \end{pmatrix}.
\end{equation*}
Note that $M=1$, which occurs when $s=0$, gives the parabolic representations, 
%that is, the traces of $\rho(x)$ and $\rho(y)$ are both $2$, 
and $\lambda=0$ gives the Alexander polynomial $\bigtriangleup(M^2)$.

 In our previous papers \cite{Jo-Kim} and \cite{Jo-Kim2}, we investigated parabolic representations of 2-bridge kmot groups.
Especially in \cite{Jo-Kim}, we converted the system of conjugation quandle equations to that of symplectic quandle equations, that is, we used the fact that the conjugation quandle ($\mathcal{D}_1$, conjugation) is isomorphic to
 the symplectic quandle $(\bc^2/{\pm 1},\langle, \rangle) $,  
 where 
% $\mathcal{D}_1$ is the set of parabolic elements of  $SL(2,\bc)$, that is, 
 $\mathcal{D}_1=\{A\in SL(2,\bc)\mid tr(A)= 2 \}$
and $\langle , \rangle$ is a symplectic form on $\bc^2$ defined by  
$\langle x,y \rangle=x_1y_2-x_2y_1$ for $ x=(x_1, x_2),\, y=(y_1,
  y_2).$ 
 
% In this paper, we study $SL(2,\bc)$-representations of 2-bridge kmot groups by extending such a symplectic quandle method.
In this paper, we extend the symplectic quandle method to study the general $SL(2,\mathbb{C})$-representations of the 2-bridge kmot groups.
%We have discovered a structure in the conjugation quandle of $\mathcal{D}_M$ that exhibits similarities to the symplectic quandle structure. As a result, we have named this structure the generalized symplectic quandle structure.
We introduce a `generalized symplectic quandle structure' corresponding to ($\mathcal{D}_M$, conjugation) for each $M\in\mathbb C\setminus \{0,1,-1\}$, where $\mathcal{D}_M=\{A\in SL(2,\bc)\mid tr(A)= M+M^{-1} \}.$ We then
  convert the system of conjugation quandle equations to that of generalized symplectic quandle equations. 
 %This time, we use the Schubert diagrams  differently from the case of parabolic representations.
As a result, we obtain a simpler expression for the 2-variable Riley polynomial and several recursive formulas for Riley polynomials and Alexander polynomials. The 2-variable Riley polynomials are simplest when they are expressed by the variables $\lambda$ and $\tilde\lambda=\lambda+s$. Note that $\lambda=2-tr(\rho(x)\rho(y)^{-1})$ and $\tilde\lambda=tr(\rho(x)\rho(y))-2$.

As an application, we can compute the A-polynomials of 2-bridge knots from our recursive formulas very efficiently.  Using Mathematica, one can obtain in a few minutes each of the  A-polynomials of 2-bridge knots with up to  12 crossings, which include many unknown A-polynomials. 
%For the double twist knots, 
%we derive an explicit recursive Resultant formula for the A-polynomials  using our $(\lambda, \tilde{\lambda})$-expressions for Riley polynomials, which is  different from that given by Petersen in \cite{Petersen}.

\section{Conjugation quandle structure on $SL(2,\bc)$} 
In this section we show how the conjugation quandle structure on $SL(2,\bc)$ can be isomorphically changed to a `generalized symplectic quandle structure' corresponding to ($\mathcal{D}_M$, conjugation). If we use the Wirtinger presentation of a knot group, an $SL(2,\bc)$ representation can be viewed as a quandle coloring, i.e., a quandle homomorphism of the knot quandle to the conjugation quandle of $SL(2,\bc)$. The actual computation in the conjugation quandle is quite difficult but one can replace the conjugation quandle by an isomorphic one which is simpler to compute and handle. Such an idea was employed in the parabolic representation case in our previous work using symplectic quandle, and we obtain a new recursive formula of the 1-variable Riley polynomial identifying it simply as a (generalized) Chebyshev polynomial. For the general representation case, computation gets more complicated but still we can do the same thing by introducing a generalized symplectic quandle structure corresponding to the conjugation quandle.

 First let us briefly review the isomorphism between  the conjugation quandle  structure on the set of parabolic elements of $ SL(2,\bc)$ and the symplectic quandle structure on $\bc^2/{\pm 1}$. We denote the set of parabolic elements of  $SL(2,\bc)$ by $\mathcal{D}_1$, that is,
$\mathcal{D}_1:=\{A\in SL(2,\bc)\mid tr(A)= 2 \}.$

%\subsection{Symplectic quandle structure on $\bc^2$}%------------------------
 Let $\langle , \rangle$ be a symplectic form on $\bc^2$ defined by  
$$\langle x,y \rangle=\begin{vmatrix}
x_1 & y_1\\
x_2 & y_2 
 \end{vmatrix}=x_1y_2-x_2y_1$$
for $ x=\begin{pmatrix}
x_1\\
  x_2
 \end{pmatrix},\, y=\begin{pmatrix}
y_1\\
  y_2
  \end{pmatrix}.$ 
Then $(\bc^2, \langle , \rangle)$ is a symplectic quandle with the quandle operation 
$$x\rhd y=x+\langle x,y \rangle y.$$
Note that if we denote the row vector $(-v,u)$ 
by  $\hat{x}$ for each colume vector $x=\begin{pmatrix}
u\\
v 
 \end{pmatrix}$, 
then $\hat{x}y=\langle x,y\rangle$ and
$$x\rhd y=x+\langle x,y \rangle y=x+(\hat{x}y)y.$$

The symplectic quandle structure on $\bc^2$ induces a quandle structure on the space of orbits of the action of multiplicative group $\{1,-1\}$ on $\bc^2$ by scalar multilication, because negating $x$ negates $x\rhd y$ and $x\rhd^{-1} y$, while  negating $y$ leaves them unchanged. We will denote the orbit space by $\mathfrak{C}=\bc^2/{\pm 1}$ and 
call this quandle $(\mathfrak{C},\langle , \rangle)$ a ($reduced$) $symplectic$ $quandle$.

%Note that there is a $\mathfrak{C}$-coloring on $\mathcal{R}(K)$ if and only if there is a $\bc^2$-coloring up to sign on $\mathcal{R}(K)$.

%\subsection{Set of parabolic elements in $SL(2,\bc)$}\label{parabolic quandle}%---------------------------
%We denote the set of parabolic elements of  $SL(2,\bc)$ by $\mathcal{D}_1$, that is, $\mathcal{D}_1=\{A\in SL(2,\bc)\mid tr(A)= 2 \}.$

From the fact that every parabolic element, which is not the identity, is conjugate to the particular element
 $\begin{pmatrix}
1 & 1\\
  0 & 1 
 \end{pmatrix}$, we obtain the following identities   
 \begin{equation*}
\begin{split}
A&=\begin{pmatrix}
a_{11}& a_{12}\\
a_{21} & a_{22} 
 \end{pmatrix}
\begin{pmatrix}
1 & 1\\
  0 & 1 
 \end{pmatrix}
\begin{pmatrix}
a_{11}& a_{12}\\
a_{21} & a_{22} 
 \end{pmatrix}^{-1}
=\begin{pmatrix}
1-a_{11}a_{21}& a_{11}^2\\
-a_{21}^2 & 1+a_{11}a_{21} 
 \end{pmatrix}\\
&=I+\begin{pmatrix}
a_{11}\\
a_{21}  
 \end{pmatrix}
(-a_{21},a_{11})=I+\begin{pmatrix}
-a_{11}\\
-a_{21}  
 \end{pmatrix}
(a_{21},-a_{11})\\
%&=I+a\hat{a}=I+(-a)\hat{(-a)}
\end{split}
\end{equation*}
and 
 \begin{equation*}
\begin{split}
A^{-1}&=I-\begin{pmatrix}
a_{11}\\
a_{21}  
 \end{pmatrix}
(-a_{21},a_{11})=I-\begin{pmatrix}
-a_{11}\\
-a_{21}  
 \end{pmatrix}
(a_{21},-a_{11})\\
%&=I-a\hat{a}=I-(-a)\hat{(-a)}.
\end{split}
\end{equation*}
This gives a bijection $T$ from  $\mathcal{D}_1$ to  $\mathfrak{C}$ such that 
$$T(A)=\left[\begin{array}{c}
a_{11}\\
a_{21}  
 \end{array}\right]:=[a]$$ 
and this map sends $A^{-1}\in \mathcal{D}_1$ to  $[ia]\in \mathfrak{C} $, that is, $T(A^{-1})=[ia]$.  Therefore we have 
$$A=I+a\hat{a},\quad A^{-1}=I-a\hat{a}=I+(ia)(\widehat{ia}).$$

The following proposition shows that $T$ defines an isomorphism between the conjugation quandle $(\mathcal{D}_1, \text{conjugation})$ and the symplectic quandle $(\mathfrak{C},\langle, \rangle) $.
\begin{proposition}[\cite{Jo-Kim}]\label{lambda-aq}
	If $T(A)=[a]$ and  $T(B)=[b]$ then 
\begin{enumerate}
\item [\rm (i)] 
$T(B^{-1}AB)=[a+\langle a,b \rangle b]\in \mathfrak{C}$
\item [\rm (ii)] 
	$T(BAB^{-1})=[a-\langle a,b \rangle b]\in \mathfrak{C}$
\item [\rm (iii)] $a\rhd b=B^{-1}a$ 
\item [\rm (iv)] $tr(AB)-2=-\langle a,b\rangle^2=\langle a,b\rangle\langle b,a\rangle $
\end{enumerate}
\end{proposition}
%\subsection{\textcolor{blue}{Some properties of $\mathcal{D}_M=\{A\in SL(2,\bc)\mid tr(A)= M+M^{-1} \}$ for $M\neq 0, \pm 1$}}%---------------------------
Now we investigate the set of elemtents of $SL(2,\bc)$ with the trace $M+M^{-1}$, $\mathcal{D}_M=\{A\in SL(2,\bc)\mid tr(A)= M+M^{-1} \}$ for $M\neq 0, \pm 1$.
The following useful lemma is well-known and its proof can be found   in Lemma 7 of \cite{Riley3}.
\begin{lemma} [\cite{Riley3}]\label{Riley-lemma-3}
Let $A, B \in \mathcal{D}_M$ and $AB\neq BA$. Then there exist $U\in SL(2,\mathbb C)$ such that
$$UAU^{-1}= \begin{pmatrix}
M & 1\\
  0& \frac{1}{M} 
 \end{pmatrix},\,\, UBU^{-1}= \begin{pmatrix}
M & 0\\
  \lambda& \frac{1}{M} 
 \end{pmatrix},$$
where $\lambda=2-tr(AB^{-1})$.
\end{lemma}
Let $z:=M-\frac{1}{M}$. Then we see that  $\begin{pmatrix}
M & 0\\
  0 & \frac{1}{M} 
 \end{pmatrix}
=
\begin{pmatrix}
1 & \frac{1}{z}\\
  0 & 1 
 \end{pmatrix}
\begin{pmatrix}
M & 1\\
  0 & \frac{1}{M} 
 \end{pmatrix}
\begin{pmatrix}
1 & -\frac{1}{z}\\
  0 & 1 
 \end{pmatrix}$ and 
   every  element in $\mathcal{D}_M$ is conjugate to the particular element
 $\begin{pmatrix}
M & 1\\
  0 & \frac{1}{M} 
 \end{pmatrix}$. And we get the following identities.   
 \begin{equation*}
\begin{split}
A&=\begin{pmatrix}
a_{11}& a_{12}\\
a_{21} & a_{22} 
 \end{pmatrix}
\begin{pmatrix}
M & 1\\
  0 & \frac{1}{M} 
 \end{pmatrix}
\begin{pmatrix}
a_{11}& a_{12}\\
a_{21} & a_{22} 
 \end{pmatrix}^{-1}
=MI+\begin{pmatrix}
a_{11}-za_{12}\\
a_{21}-za_{22}  
 \end{pmatrix}
(-a_{21},a_{11})\\
&=\frac{1}{M} I+\begin{pmatrix}
a_{11}\\
a_{21}  
 \end{pmatrix}
(-a_{21}+za_{22},a_{11}-za_{12}),
\end{split}
\end{equation*}
 If we denote the first and the second colume vector of the matrix $\begin{pmatrix}
a_{11}& a_{12}\\
a_{21} & a_{22} 
 \end{pmatrix}$ by $a_1$ and $a_2$ respectively and denote the vector $a_1-za_2$ by $a_3$, then 
 $A$ and $A^{-1}$ can be expressed as
\begin{equation}\label{quandle-1}
\begin{split}
A&=MI+a_3\hat{a_1}=\frac{1}{M}I+a_1\hat{a_3},\\
A^{-1}&=MI-a_1\hat{a_3}=\frac{1}{M}I-a_3\hat{a_1}.
\end{split}
\end{equation}
We will denote the matrix $(a_1,a_3)$ by $\tilde{a}$ from now on, and this will play the role of a $\bc^2$ vector $a$ of the parabolic case.  %$\det\tilde{a}=-z$.
Note that the matrix $\tilde{a}$ is an element of the set 
 $$\tilde{\mathcal{D}}_M:=\{A\in GL(2,\bc)\mid det(A)=- M+M^{-1} \}$$
for each $a\in SL(2,\mathbb C)$, since
 $$\det\tilde{a}=\langle a_1, a_3 \rangle=\langle a_1, a_1-z a_2 \rangle=-z\langle a_1, a_2 \rangle=-z=-M+\frac{1}{M}.$$

%\subsection{Generalized symplectic quandle structure corresponding to ($\mathcal{D}_M$, conjugation)}%---------------------------
%We have seen in section \ref{parabolic quandle} that  the conjugation qundle ($\mathcal{D}_1$, conjugation) can be considered as a symplectic quandle. In this section 

Now we show that there is a similar structure on ($\mathcal{D}_M$,conjugation) for each $M\in\mathbb C\setminus \{0,1,-1\}$, which we will call {\it a generalized symplectic quandle}. Firstly, we compute the isotropy subgroup $S_M$ of $\begin{pmatrix}
M & 1\\
  0 & \frac{1}{M} 
 \end{pmatrix}$ in $SL(2,\bc)$ under the action of conjugation.
\begin{proposition}
Let	 $S_M=\{	\begin{pmatrix}
t & \frac{z'}{z}\\
  0 & \frac{1}{t} 
 \end{pmatrix}
\mid z=M-\frac{1}{M}, z'=t-\frac{1}{t}\}$. Then $$X\in S_M\,\, \text{ if and only if} \,\,\,X\begin{pmatrix}
M & 1\\
  0 & \frac{1}{M} 
 \end{pmatrix}X^{-1}=\begin{pmatrix}
M & 1\\
  0 & \frac{1}{M} 
 \end{pmatrix}.$$
\end{proposition}
\begin{proof}
It is obvious that $X\begin{pmatrix}
M & 1\\
  0 & \frac{1}{M} 
 \end{pmatrix}X^{-1}=\begin{pmatrix}
M & 1\\
  0 & \frac{1}{M} 
 \end{pmatrix}$ if $X\in S_M$. 
Let $X=\begin{pmatrix}
a & b\\
  c & d 
 \end{pmatrix}$
and suppose that $X\begin{pmatrix}
M & 1\\
  0 & \frac{1}{M} 
 \end{pmatrix}X^{-1}=\begin{pmatrix}
M & 1\\
  0 & \frac{1}{M} 
 \end{pmatrix}$. Then from
$$\frac{1}{M}I+X_1\hat{X_3}=X\begin{pmatrix}
M & 1\\
  0 & \frac{1}{M} 
 \end{pmatrix}X^{-1}=\begin{pmatrix}
M & 1\\
  0 & \frac{1}{M} 
 \end{pmatrix}=\frac{1}{M}I+
\begin{pmatrix}
z & 1\\
  0 & 0 
 \end{pmatrix}$$
we have
$$\begin{pmatrix}
z & 1\\
  0 & 0 
 \end{pmatrix}=X_1\hat{X_3}=
\begin{pmatrix}
a\\
c 
 \end{pmatrix}
(-c+zd,a-zb)=z
\begin{pmatrix}
ad & -ab\\
  cd & -bc 
 \end{pmatrix}
+\begin{pmatrix}
-ac & a^2\\
  -c^2 & ac 
 \end{pmatrix},
$$
which implies that $$c=0,\,\, d=\frac{1}{a}, \,\,zb=a-\frac{1}{a}.$$
Hence 
$X=\begin{pmatrix}
a & \frac{zb}{z}\\
  0 & \frac{1}{a} 
 \end{pmatrix}\in S_M$.
\end{proof}
 Now we define an equivalence relation on  $SL(2,\bc)$  and  another  equivalence relation on 
 $\tilde{\mathcal{D}}_M$ 
respectively, and then look at the correspondence between these two equivalence relations. 
\begin{definition} Let $a, b \in SL(2,\bc)$ and  $q(t)=\begin{pmatrix}
t & \frac{t-t^{-1}}{z}\\
  0 & t^{-1} 
 \end{pmatrix}, \,\,p(t)=\begin{pmatrix}
t & 0\\
  0 &  t^{-1}
 \end{pmatrix}$.
Two equivalence relations and the {\it quandle matrix $[\tilde{a},\tilde{b}]$ for $a$ and $b$} are defined as follows.
\begin{enumerate}
\item [\rm (i)]  $a\sim b$ if  $aq(t)=b$ for some $t\in \mathbb R$, or equivalently, $a^{-1}b\in S_M$,
\item [\rm (ii)] $\tilde{a}\sim \tilde{b}$  if $\tilde{a}p(t)=\tilde{b} $ for some $t\in \mathbb R$,
\item [\rm (iii)] $[\tilde{a},\tilde{b}]=\begin{pmatrix}
\langle a_1, b_3 \rangle& 0\\
  0 & \langle a_3, b_1 \rangle
 \end{pmatrix}$%, where $\langle a_i. b_j \rangle:=\det (a_i,b_j)$.
\end{enumerate}
\end{definition}
\begin{remark}\label{a-a}
Note that 
$[\tilde{a},\tilde{a}]=
\begin{pmatrix}
-z & 0\\
  0 & z
 \end{pmatrix}$ for any $a \in SL(2,\bc)$ since  $\langle a_1, a_3 \rangle=\det\tilde{a}=-z,$ and it follows easily from the definition that  
$$[\tilde{a},\tilde{b}]=[y\tilde{a},y\tilde{b}]=[\widetilde{ya},\widetilde{yb}]$$ for any $y \in SL(2,\bc)$.
\item
%Note that if $M=1$, then we have
%$[\tilde{a},\tilde{b}]
%\begin{pmatrix}\langle a_1, b_3 \rangle& 0\\  0 & \langle a_3, b_1 \rangle \end{pmatrix}=\begin{pmatrix}\langle a_1, b_1\rangle& 0\\  0 & \langle a_1, b_1 \rangle \end{pmatrix}
%=\langle a_1, b_1 \rangle Id$, which is the
\end{remark}
From now on, we will use  the notation $m$ and $e$ respectively for the particular matrices 
$\begin{pmatrix}
M & 1\\
  0 & \frac{1}{M} 
 \end{pmatrix}$ and 
$\begin{pmatrix}
M & 0\\
  0 & \frac{1}{M} 
 \end{pmatrix}$. 
\begin{lemma}
$a\sim b$ if and only if $\tilde{a}\sim \tilde{b}$.
\end{lemma}
\begin{proof}
$a\sim b$ if and only if there exists $t\in\mathbb R$ such that 
$$b=aq(t)=(a_1,a_2)\begin{pmatrix}
t & \frac{z'}{z}\\
  0 & \frac{1}{t} 
 \end{pmatrix}=(ta_1,\frac{z'}{z}a_1+\frac{1}{t} a_2),\,\,z'=t-t^{-1}.$$
Since
\begin{equation*}
\begin{split}
b_3&=ta_1-z(\frac{z'}{z}a_1+\frac{1}{t} a_2)=ta_1-z'a_1-\frac{z}{t}a_2\\
&=(t-z')a_1-\frac{z}{t}a_2=\frac{1}{t}a_1-\frac{z}{t}a_2=\frac{1}{t}(a_1-za_2)=\frac{1}{t}a_3,
\end{split}
\end{equation*}
 we have $$\tilde{b}=(b_1,b_3)=(ta_1, \frac{1}{t}a_3)=\tilde{a}p(t),$$ which completes the proof.
\end{proof}
\begin{corollary}\label{quandle-equation}
Let $a, b \in SL(2,\bc)$ and  $x=ama^{-1},\,\, y=bmb^{-1}$. Then 
$$x=y \,\,\, \Longleftrightarrow\,\,\, a \sim b\,\,\, \Longleftrightarrow \,\,\,
\tilde{a}\sim \tilde{b}.$$
\end{corollary}

From (\ref{quandle-1}) and  the above corollary, 
we have a bijection $T$ from 
$\mathcal{D}_M$ to $\tilde{\mathcal{D}}_M/\sim$ such that 
$$T(x)=T(ama^{-1})=[\tilde{a}].$$
%where $$\tilde{\mathcal{D}}_M=\{A\in GL(2,\bc)\mid det(A)=- M+M^{-1} \}.$$ 
So we can think  of $x$ as an equivalence class of $\tilde{a}$ and this induces 
a `generalized symplectic quandle structure' corresponding to ($\mathcal{D}_M$, conjugation), which is similar to the symplectic quandle structure corresponding to ($\mathcal{D}_1$, conjugation). 

\begin{lemma}\label{quandle-comp-1}
Let $\epsilon$ be either $1$ or $-1$ and  $w=x\rhd^{\epsilon} y:=y^{-\epsilon}xy^\epsilon$.  If $x=ama^{-1}$, $y=bmb^{-1}$, and  $w=cmc^{-1}$, then  $c\sim y^{-\epsilon}a$ and $\tilde{c}\sim y^{-\epsilon}\tilde{a}$.
\end{lemma}
\begin{proof} Since 
\begin{equation*}
\begin{split}
w&=y^{-\epsilon}xy^\epsilon=(bmb^{-1})^{-\epsilon}ama^{-1}(bmb^{-1})^{\epsilon}=(bm^{-\epsilon}b^{-1}a)m(bm^{-\epsilon}b^{-1}a)^{-1},\\
\end{split}
\end{equation*}
we have 
$c\sim bm^{-\epsilon}b^{-1}a=y^{-\epsilon}a$ 
and thus 
$\tilde{c}\sim \widetilde{y^{-\epsilon}a}\sim y^{-\epsilon}\tilde{a}$.
\end{proof}
Now if we denote the class $[\tilde{c}]$ by $\tilde{a}\rhd \tilde{b}$ in the case when $x=ama^{-1}$, $y=bmb^{-1}$, $w=cmc^{-1}$,  and $w=x\rhd y$, then $(\tilde{\mathcal{D}}_M/\sim, \rhd)$ becomes a quandle.  
\begin{proposition}\label{quandle-axiom}
$(\tilde{\mathcal{D}}_M/\sim, \rhd)$ is a quandle.
\end{proposition}
\begin{proof}
Since
\begin{equation*}
\begin{split}
\centerdot\quad &\tilde{a}\rhd \tilde{a} \sim \widetilde{am^{-1}a^{-1}a}\sim\widetilde{am^{-1}}\sim \tilde{a}\\
\centerdot\quad &(\tilde{a}\rhd \tilde{b}) \rhd \tilde{c}\sim \widetilde{bm^{-1}b^{-1}a}\rhd \tilde{c} \sim \widetilde{cm^{-1}c^{-1}bm^{-1}b^{-1}a}\\
\centerdot\quad &(\tilde{a}\rhd \tilde{c}) \rhd (\tilde{b}\rhd \tilde{c})\sim \widetilde{cm^{-1}c^{-1}a}\rhd \widetilde{cm^{-1}c^{-1}b}   \sim \widetilde{cm^{-1}c^{-1}bm^{-1}b^{-1}a}
\end{split}
\end{equation*}
by Lemma \ref{quandle-comp-1}, 
(i) $\tilde{a}\rhd \tilde{a} \sim \tilde{a}$ and (ii) $(\tilde{a}\rhd \tilde{b}) \rhd \tilde{c}\sim (\tilde{a}\rhd \tilde{c}) \rhd (\tilde{b}\rhd \tilde{c})$ hold. If we let $c=bmb^{-1}a$, then (iii) $\tilde{c}\rhd \tilde{b}\sim  \tilde{a}$ also holds because
$$\tilde{c}\rhd \tilde{b} \sim \widetilde{bm^{-1}b^{-1}c}\sim \widetilde{bm^{-1}b^{-1}bmb^{-1}a}\sim  \tilde{a},  $$
and such $\tilde{c}$ is unique, that is, if there is another $\tilde{c'}$ such that  $\tilde{c'}\rhd \tilde{b}\sim  \tilde{a}$ then  $\tilde{c'}\sim \tilde{c}$. 
Therefore $(\tilde{\mathcal{D}}_M/\sim, \rhd)$ satisfies the quandle axioms as desired.
\end{proof}
Up to now, we have defined a quandle structure on $\tilde{\mathcal{D}}_M/\sim$ which is isomorphic to the conjugation quandle structure on $\mathcal{D}_M$.
The following proposition  shows that the quandle  $(\tilde{\mathcal{D}}_M/\sim, \rhd)$ is very similar to the symplectic quandle  $(\mathfrak{C},\langle, \rangle) $ 
so that  we will call $(\tilde{\mathcal{D}}_M/\sim, \rhd)$ and $[ \,,\, ]$    {\it a generalized symplectic quandle} and {\it a generalized symplectic form} respectively, and
denote it by  $(\tilde{\mathcal{D}}_M/\sim, [ \,,\, ])$. Indeed, $tr[\tilde{a},\tilde{b}]$ is a symplectic form. 
\begin{proposition}\label{quandle-method}
Let $\epsilon$ be either $1$ or $-1$ and  $w=x\rhd^{\epsilon} y=y^{-\epsilon}xy^\epsilon$.  If $x=ama^{-1}$ and $y=bmb^{-1}$, then 
$$w=cmc^{-1} \,\,\, \text{with}\,\,\,\tilde{c}\sim \tilde{a}e^{\epsilon}+\epsilon \tilde{b}[\tilde{a}, \tilde{b}],$$ 
\end{proposition}
\begin{proof}
Suppose that $\epsilon=1$. Then by lemma \ref{quandle-comp-1} and (\ref{quandle-1}),
\begin{equation*}
\begin{split}
\tilde{c}&\sim y^{-1}\tilde{a}=(MI-b_1\hat{b_3})(a_1,a_3)=(Ma_1+\langle a_1,b_3\rangle b_1, Ma_3+\langle a_3,b_3\rangle b_1)\\
\end{split}
\end{equation*}
and
\begin{equation*}
\begin{split}
\tilde{c}&\sim y^{-1}\tilde{a}=(\frac{1}{M}-b_3\hat{a_1})(a_1,a_3)=(\frac{1}{M}a_1+\langle a_1,b_1\rangle b_3, \frac{1}{M}a_3+\langle a_3,b_1\rangle b_3).\\
\end{split}
\end{equation*}
Hence we have
\begin{equation*}
\begin{split}
\tilde{c}\sim y^{-1}\tilde{a}&=(Ma_1+\langle a_1,b_3\rangle b_1,\frac{1}{M}a_3+\langle a_3,b_1\rangle b_3)\\
&=(a_1,a_3)\begin{pmatrix}
M & 0\\
 0& \frac{1}{M} 
 \end{pmatrix}+(b_1,b_3)\begin{pmatrix}
\langle a_1, b_3 \rangle& 0\\
  0 & \langle a_3, b_1 \rangle
 \end{pmatrix}\\
&=\tilde{a}\begin{pmatrix}
M & 0\\
 0& \frac{1}{M} 
 \end{pmatrix}+\tilde{b}\begin{pmatrix}
\langle a_1, b_3 \rangle& 0\\
  0 & \langle a_3, b_1 \rangle
 \end{pmatrix}\\
&=\tilde{a}e+\tilde{b}[\tilde{a}, \tilde{b}].
\end{split}
\end{equation*}
Similarly, we have $\tilde{c}\sim y\tilde{a}=\tilde{a}e^{-1}-\tilde{b}[\tilde{a}, \tilde{b}]$ when $\epsilon= -1$.
\end{proof}
\begin{remark}
Note that if $M=1$, then we recover the old symplectic quandle equation as follows:
\begin{equation*}
\begin{split}
(c_1,c_1)=\tilde{c}&\sim
\tilde{a}e^{\epsilon}+\epsilon \tilde{b}[\tilde{a}, \tilde{b}]\\
&=(a_1,a_1)+\epsilon (b_1,b_1) 
\begin{pmatrix}
\langle a_1, b_1 \rangle& 0\\
  0 & \langle a_1, b_1 \rangle
 \end{pmatrix}\\
&=(a_1+\epsilon\langle a_1, b_1 \rangle b_1,\,\, a_1+\epsilon\langle a_1, b_1 \rangle b_1)
 \end{split}
\end{equation*}
and $c_1 \sim a_1 +\epsilon\langle a_1, b_1 \rangle b_1$ is the symplectic quandle equation corresponding to $w=x\rhd^{\epsilon} y$ when $x$ and $y$ are parabolic elements in $SL(2,\mathbb C)$ (see \cite{Jo-Kim}).
\end{remark}

Now we denote  $z^2+\lambda$ by $\tilde{\lambda}$.
Then 
$$tr(\begin{pmatrix}
M & 1\\
  0& \frac{1}{M} 
 \end{pmatrix} \begin{pmatrix}
M & 0\\
  \lambda& \frac{1}{M} 
 \end{pmatrix})=2+z^2+\lambda=2+\tilde{\lambda}.$$ 
%\jocomment{\textcolor{red}{geometric meaning?}}
and we get the following which will be of great importance later. Note that the  identity in the next lemma corresponds to the identity (iv) of Proposition \ref{lambda-aq} since $\tilde{\lambda}=\lambda$ when $M=1$.
%$[\tilde{a},\tilde{b}][\tilde{b},\tilde{a}]=\begin{pmatrix}\tilde{\lambda}  & 0\\  0& \tilde{\lambda} \end{pmatrix}$ for any  $a,b \in SL(2,\mathbb C)$.

\begin{lemma}\label{a-a-b}
Let  $A=ama^{-1}, B=bmb^{-1}$ with $a,b \in SL(2,\mathbb C)$.
Then we have
$$
 [\tilde{a},\tilde{b}][\tilde{b},\tilde{a}]=\begin{pmatrix}
\tilde{\lambda}  & 0\\
  0& \tilde{\lambda}
 \end{pmatrix},\,\, \tilde{\lambda}=tr(AB)-2
$$ if $AB\neq BA$.
\end{lemma}
\begin{proof}
By Lemma \ref{Riley-lemma-3}, we have $U\in SL(2,\mathbb C)$ such that 
$$UAU^{-1}= \begin{pmatrix}
M & 1\\
  0& \frac{1}{M} 
 \end{pmatrix},\,\, UBU^{-1}= \begin{pmatrix}
M & 0\\
  \lambda& \frac{1}{M} 
 \end{pmatrix},\,\,\lambda=2-tr(AB^{-1}).$$
 That is, 
 $$m=(Ua)m(Ua)^{-1},\,\, (Ub)m(Ub)^{-1}= \begin{pmatrix}
M & 0\\
  \lambda& \frac{1}{M} 
 \end{pmatrix}.$$
Hence 
$Ua= q(t)=\begin{pmatrix}
t & \frac{t-t^{-1}}{z}\\
  0 & t^{-1} 
 \end{pmatrix},\,\,\, Ub= \begin{pmatrix}
z & 1\\
\lambda& \frac{\lambda+1}{z}
 \end{pmatrix}q(t')$ for some nonzero $t,t'$,
since
\begin{equation}\label{b}
\begin{pmatrix}
z & 1\\
\lambda& \frac{\lambda+1}{z}
 \end{pmatrix}m\begin{pmatrix}
z & 1\\
\lambda& \frac{\lambda+1}{z}
 \end{pmatrix}^{-1}= \begin{pmatrix}
M & 0\\
  \lambda& \frac{1}{M} 
 \end{pmatrix}.
\end{equation}

 Now we have
\begin{equation*}
\begin{split}
 [\tilde{a},\tilde{b}]&= [\widetilde{Ua},\widetilde{Ub}]
 = [\begin{pmatrix}
t & \frac{1}{t}\\
  0& -\frac{z}{t} 
 \end{pmatrix},\begin{pmatrix}
zt' & 0\\
  \lambda t'& -\frac{1}{t'} 
 \end{pmatrix}]=\begin{pmatrix}
-\frac{t}{t'} & 0\\
  0& \frac{t'}{t}\tilde{\lambda}
 \end{pmatrix}=\begin{pmatrix}
-1 & 0\\
  0& \tilde{\lambda} 
 \end{pmatrix}p(\frac{t}{t'})\\
 [\tilde{b},\tilde{a}]&= [\widetilde{Ub},\widetilde{Ua}] = [\begin{pmatrix}
zt' & 0\\
  \lambda t'& -\frac{1}{t'} 
 \end{pmatrix},\begin{pmatrix}
t & \frac{1}{t}\\
  0& -\frac{z}{t} 
 \end{pmatrix}]=\begin{pmatrix}
-\frac{t'}{t}\tilde{\lambda} & 0\\
  0& \frac{t}{t'}
 \end{pmatrix}=\begin{pmatrix}
-\tilde{\lambda}  & 0\\
  0& 1
 \end{pmatrix}p(\frac{t'}{t}),\\
\end{split}
\end{equation*}
which implies $ [\tilde{a},\tilde{b}][\tilde{b},\tilde{a}]=\begin{pmatrix}
\tilde{\lambda}  & 0\\
  0& \tilde{\lambda}
 \end{pmatrix}.$
\end{proof}

\section{$\epsilon_i$-sequences and $c_k^n, d_k^n, \tilde{c}_k^n, \tilde{d}_k^n$} %===================================
In order to study $SL(2,\bc)$-representations of 2-bridge kmot groups and to obtain an explicit formula for the Riley polynomial in the next section, we need to introduce some families of polynomials in $\mathbb Z[M^{\pm 1}]$ for a given sequence $\epsilon=\{\epsilon_i\,|\, \epsilon_i=1\,\,\text{or}\,\, \epsilon_i=-1,\,\,i\in\mathbb Z\}$.
\begin{definition}\label{cd}  For a given sequence $\epsilon=\{\epsilon_i\,|\, 
\epsilon_i=1\,\,\text{or}\,\, \epsilon_i=-1,\,\,
i\in\mathbb Z\}$, we define $c_k^n, d_k^n, \tilde{c}_k^n, \tilde{d}_k^n$ as follows.
\begin{enumerate}
\item [\rm (i)]
$\displaystyle c_k^n(M,\epsilon)=\sum_{\shortstack{$\scriptstyle i_1:\text{odd} $\\$\scriptstyle  i_k\leq n$}}^{\wedge} \epsilon_{i_1} \cdots\epsilon_{i_{k}}M^{\widehat{i_1 \cdots i_{k}}}, \quad
\displaystyle d_k^n(M,\epsilon)=\sum_{\shortstack{$\scriptstyle i_1:\text{even} $\\$\scriptstyle  i_k\leq n$}}^{\wedge} \epsilon_{i_1} \cdots\epsilon_{i_{k}}M^{-\widehat{i_1 \cdots i_{k}}}$
\item [\rm (ii)]
$\displaystyle \tilde{c}_k^n(M,\epsilon)=\sum_{\shortstack{$\scriptstyle i_1:\text{odd} $\\$\scriptstyle  i_k\leq n$}}^{\wedge} \epsilon_{i_1} \cdots\epsilon_{i_{k}}M^{\widehat{i_1 \cdots i_{k}}^{alt}}, \quad
\displaystyle \tilde{d}_k^n(M,\epsilon)=\sum_{\shortstack{$\scriptstyle i_1:\text{even} $\\$\scriptstyle  i_k\leq n$}}^{\wedge} \epsilon_{i_1} \cdots\epsilon_{i_{k}}M^{-\widehat{i_1 \cdots i_{k}}^{alt}}$.
\end{enumerate}
Here $\displaystyle\sum_{\shortstack{$\scriptstyle i_1:\text{odd} $\\$\scriptstyle  i_k\leq n$}}^{\wedge} $,  $\displaystyle\sum_{\shortstack{$\scriptstyle i_1:\text{even} $\\$\scriptstyle  i_k\leq n$}}^{\wedge} $, $\widehat{i_1 \cdots i_{k}}$, and $\widehat{i_1 \cdots i_{k}}^{alt}$ mean the followings.
\begin{equation*}
\begin{split}
\displaystyle\sum_{\shortstack{$\scriptstyle i_1:\text{odd} $\\$\scriptstyle  i_k\leq n$}}^{\wedge}&  :0< i_1< i_2 < i_3 <\cdots<i_k\leq n, \,\, i_{2i+1} \text{ is odd},\,\,  i_{2i} \text{ is even}\\
\displaystyle\sum_{\shortstack{$\scriptstyle i_1:\text{even} $\\$\scriptstyle  i_k\leq n$}}^{\wedge}&  : 0<i_1< i_2 < i_3 <\cdots<i_k\leq n, \,\, i_{2i+1} \text{ is even},\,\,  i_{2i} \text{ is odd}\\
\widehat{i_1 \cdots i_{k}}&= \epsilon_{1}+\cdots+ \epsilon_{i_1-1}-( \epsilon_{i_1+1}+\cdots+ \epsilon_{i_2-1})+\cdots+(-1)^{k}( \epsilon_{i_k+1}+\cdots+\epsilon_{n})\\
\widehat{i_1 \cdots i_{k}}^{alt}&= -\epsilon_{1}+\epsilon_2-\cdots +(-1)^{n-k} \epsilon_{n}\,\,(\text{alternating sum of }\,\,\{\epsilon_1, \epsilon_2, \cdots, \epsilon_{n}\}-\{\epsilon_{i_1}, \epsilon_{i_2},\cdots, \epsilon_{i_k}\})
\end{split}
\end{equation*}
For example, $\widehat{258}= \epsilon_{1}-\epsilon_3-\epsilon_4+\epsilon_6+\epsilon_7-\epsilon_{9}-\cdots -\epsilon_{n}$ and $\widehat{2367}^{alt}= -\epsilon_{1}+\epsilon_4-\epsilon_5+\epsilon_8-\epsilon_9+\cdots +\epsilon_{n}$.
\end{definition}
For the sake of convenience,  we define $c^{n}_{-1}, \tilde{c}^{n}_{-1}, c^{n}_{0}, \tilde{c}^{n}_{0} $ and $  d^{n}_{-1}, \tilde{d}^{n}_{-1},  d^{n}_{0}, \tilde{d}^{n}_{0}, $  as follows:
$$c^{n}_{-1}:=0, \tilde{c}^{n}_{-1}:=0, c^{n}_{0}:=M^{\epsilon_{1}+\epsilon_{2}+\cdots+\epsilon_{n} },\tilde{c}^{n}_{0}:=M^{-\epsilon_{1}+\epsilon_{2}+\cdots+(-1)^n\epsilon_{n} },$$
$$d^{n}_{-1}:=0, \tilde{d}^{n}_{-1}:=0, d^{n}_{0}:=M^{-(\epsilon_{1}+\epsilon_{2}+\cdots+\epsilon_{n} )}, \tilde{d}^{n}_{0}:=M^{\epsilon_{1}-\epsilon_{2}+\cdots+(-1)^{n+1}\epsilon_{n} }.$$
These will be useful in expressing certain identities more conveniently.

Note that  for any $\square\in \{c, d, \tilde{c}, \tilde{d}\}$
 $$ \square_{2k}^{2n}(M^{-1}, \epsilon)= \square_{2k}^{2n}(M, -\epsilon)$$ and $$\square_{2k+1}^{2n}(M^{-1}, \epsilon)= -\square_{2k+1}^{2n}(M, -\epsilon).$$
 From the definitions of  $c_k^n, d_k^n, \tilde{c}_k^n, \tilde{d}_k^n$, we have the following lemmas immediately.
\begin{lemma}\label{lem-1}
\begin{enumerate}
%\item [\rm (i)]$c_{2k+1}^{2n}=c_{2k+1}^{2n-1}M^{-\epsilon_{2n}}$,  $d_{2k}^{2n+1}=d_{2k}^{2n}M^{-\epsilon_{2n+1}}$
\item [\rm (i)]
$c_{2k+1}^{2n}=c_{2k+1}^{2n-1}M^{-\epsilon_{2n}}$, $c_{2k}^{2n+1}=c_{2k}^{2n}M^{\epsilon_{2n+1}}$
\item [\rm (ii)]$c_{2k}^{2n}=c_{2k}^{2n-2}M^{\epsilon_{2n-1}+\epsilon_{2n}}+\epsilon_{2n} c_{2k-1}^{2n-1} $
\item [\rm (iii)]
$c_{2k+1}^{2n+1}=c_{2k+1}^{2n-1}M^{-\epsilon_{2n}-\epsilon_{2n+1}}+\epsilon_{2n+1} c_{2k}^{2n} $
\item [\rm (iv)]
$d_{2k}^{2n}=d_{2k}^{2n-1}M^{-\epsilon_{2n}}$, $d_{2k+1}^{2n+1}=d_{2k+1}^{2n}M^{\epsilon_{2n+1}}$
\item [\rm (v)]$d_{2k+1}^{2n}=d_{2k+1}^{2n-2}M^{\epsilon_{2n-1}+\epsilon_{2n}}+\epsilon_{2n} d_{2k}^{2n-1} $
\item [\rm (vi)]
$d_{2k}^{2n+1}=d_{2k}^{2n-1}M^{-\epsilon_{2n}-\epsilon_{2n+1}}+\epsilon_{2n+1} d_{2k-1}^{2n} $
\end{enumerate}
\end{lemma}

\begin{lemma}\label{symmetric}
If $\epsilon=(\epsilon_1, \epsilon_2,\cdots, \epsilon_{2n})$ is symmetric, then
\begin{enumerate}
    \item  [\rm (i)] $ c_{2k+1}^{2n}(M, \epsilon)= d_{2k+1}^{2n}(M, \epsilon)$
 \item  [\rm (ii)] $ \tilde{c}_{2k+1}^{2n}(M, \epsilon)= \tilde{d}_{2k+1}^{2n}(M^{-1}, \epsilon)$
    \item  [\rm (iii)] $\tilde{c}_{2k}^{2n}(M^{-1}, \epsilon)=\tilde{c}_{2k}^{2n}(M, \epsilon)=\tilde{c}_{2k}^{2n}(M, -\epsilon)$
    %\item $\tilde{c}_{2k+1}^{2n}(M^{-1}, \epsilon)=-\tilde{c}_{2k+1}^{2n}(M, -\epsilon)$
\end{enumerate}
\end{lemma}
Note that $c_{2k}^{2n}(M^{-1}, \epsilon)\neq c_{2k}^{2n}(M, \epsilon)$ even though $\epsilon$ is symmetric. 

\begin{lemma}\label{lem-2}
%$c_k^n$ and $d_k^n$ satisfy the followings.
\begin{enumerate}
\item [\rm (i)]
$c_{2k}^{2l-2}M^{\epsilon_{2l-1}+\epsilon_{2l}}  +\epsilon_{2l-1}\epsilon_{2l} c_{2k-2}^{2l-2} +\epsilon_{2l} c_{2k-1}^{2l-2}M^{-\epsilon_{2l-1}}=c_{2k}^{2l}$
\item [\rm (ii)]
$c_{2k+1}^{2l-2}M^{-\epsilon_{2l-1}-\epsilon_{2l}}+\epsilon_{2l-1}c_{2k}^{2l-2}M^{-\epsilon_{2l}}=c_{2k+1}^{2l}$
\item [\rm (iii)]
$d_{2k+1}^{2l-2}M^{\epsilon_{2l-1}+\epsilon_{2l}}+\epsilon_{2l}d_{2k}^{2l-2}M^{-\epsilon_{2l-1}}+ \epsilon_{2l-1}\epsilon_{2l}d_{2k-1}^{2l-2}=d_{2k+1}^{2l}$
%$d_{2k+1}^{2l-2}+\epsilon_{2l}d_{2k}^{2l-2}+ \epsilon_{2l-1}\epsilon_{2l}d_{2k-1}^{2l-2}=d_{2k+1}^{2l}$
\item [\rm (iv)]
$d_{2k}^{2l-2}M^{-\epsilon_{2l-1}-\epsilon_{2l}}+\epsilon_{2l-1}d_{2k-1}^{2l-2}M^{-\epsilon_{2l}}=d_{2k}^{2l}$
\item [\rm (v)]
$\tilde{c}_{2k}^{2l-2}M^{-\epsilon_{2l-1}+\epsilon_{2l}} +\epsilon_{2l}M^{\epsilon_{2l-1}} \tilde{c}_{2k-1}^{2l-2} +\epsilon_{2l-1}\epsilon_{2l} \tilde{c}_{2k-2}^{2l-2} =\tilde{c}_{2k}^{2l}$
\item [\rm (vi)]
$\tilde{c}_{2k+1}^{2l-2}M^{\epsilon_{2l-1}-\epsilon_{2l}}+\epsilon_{2l-1}\tilde{c}_{2k}^{2l-2}M^{-\epsilon_{2l}}=\tilde{c}_{2k+1}^{2l}$
\item [\rm (vii)]
$\tilde{d}_{2k+1}^{2l-2}M^{-\epsilon_{2l-1}+\epsilon_{2l}}+\epsilon_{2l}\tilde{d}_{2k}^{2l-2}M^{\epsilon_{2l-1}}+ \epsilon_{2l-1}\epsilon_{2l}\tilde{d}_{2k-1}^{2l-2}=\tilde{d}_{2k+1}^{2l}$
\item [\rm (viii)] $\tilde{d}_{2k}^{2l-2}M^{\epsilon_{2l-1}-\epsilon_{2l}}+\epsilon_{2l-1}\tilde{d}_{2k-1}^{2l-2}M^{-\epsilon_{2l}}=\tilde{d}_{2k}^{2l}$
\end{enumerate}
\end{lemma}
\begin{proof}The identities (i) and (ii) are proved using Lemma \ref{lem-1} as follows:
\begin{equation*}
\begin{split}
\rm (i) \quad c_{2k}^{2l}&=c_{2k}^{2l-2}M^{\epsilon_{2l-1}+\epsilon_{2l}}+\epsilon_{2l} c_{2k-1}^{2l-1}=c_{2k}^{2l-2}M^{\epsilon_{2l-1}+\epsilon_{2l}} +\epsilon_{2l}(c_{2k-1}^{2l-3}M^{-\epsilon_{2l-1}-\epsilon_{2l-2}}+\epsilon_{2l-1} c_{2k-2}^{2l-2})\\
&=c_{2k}^{2l-2}M^{\epsilon_{2l-1}+\epsilon_{2l}}  +\epsilon_{2l-1}\epsilon_{2l} c_{2k-2}^{2l-2} +\epsilon_{2l} c_{2k-1}^{2l-3}M^{-\epsilon_{2l-1}-\epsilon_{2l-2}}\\
&=c_{2k}^{2l-2} M^{\epsilon_{2l-1}+\epsilon_{2l}} +\epsilon_{2l-1}\epsilon_{2l} c_{2k-2}^{2l-2} +\epsilon_{2l} c_{2k-1}^{2l-2}M^{-\epsilon_{2l-1}}\\
&\\
\rm(ii)\quad c_{2k+1}^{2l}&=c_{2k+1}^{2l-1}M^{-\epsilon_{2l}}=(c_{2k+1}^{2l-3}M^{-\epsilon_{2l-2}-\epsilon_{2l-1}}+\epsilon_{2l-1}c_{2k}^{2l-2})M^{-\epsilon_{2l}}\\
&=c_{2k+1}^{2l-2}M^{-\epsilon_{2l-1}-\epsilon_{2l}}+\epsilon_{2l-1}c_{2k}^{2l-2}M^{-\epsilon_{2l}}
\end{split}
\end{equation*}
The identities from (iii) to (viii) can be similarly proved, and we omit their proofs here.
\end{proof}

The following proposition tells us the relationship between $c_k^n, d_k^n$ and $\tilde{c}_k^n, \tilde{d}_k^n$. 
\begin{proposition}\label{prop-1}
Let $\tilde{\lambda}=\lambda+z^2,\,\, z=M-M^{-1}$ and $\sigma^-$
be a shifting map on $\{\epsilon_i\,|\, i\in\mathbb Z\}$ defined by  $\sigma^-(\epsilon_i)=\epsilon_{i-1}.$ Then 
\begin{enumerate}
\item [\rm (i)]
$\displaystyle\sum_{k=0}^{l-1} \tilde{c}_{2k+1}^{2l} \tilde{\lambda}^k=\displaystyle\sum_{k=0}^{l-1} c_{2k+1}^{2l} \lambda^k$ %(if $\{\epsilon_i\}$ is symmetric?)
\item [\rm (ii)]
$\displaystyle\sum_{k=0}^{l-1} \tilde{d}_{2k+1}^{2l} \tilde{\lambda}^k=\displaystyle\sum_{k=0}^{l-1} d_{2k+1}^{2l}(M^{-1},{\epsilon}) \lambda^k$ 
\item [\rm (iii)]
$\displaystyle\sum_{k=0}^l \tilde{c}_{2k}^{2l} \tilde{\lambda}^k=\displaystyle\sum_{k=0}^l c_{2k}^{2l} \lambda^k-z\displaystyle\sum_{k=0}^{l-1} c_{2k+1}^{2l} \lambda^k=\displaystyle\sum_{k=0}^l (c_{2k}^{2l} -zc_{2k+1}^{2l}) \lambda^k$
\item [\rm (iv)]
$\displaystyle\sum_{k=0}^{l-1} d_{2k}^{2l} \lambda^k
=M^{-\epsilon_1-\epsilon_{2l}}
\displaystyle\sum_{k=0}^{l-1} (\tilde{c}_{2k}^{2l-2}
-z\tilde{c}_{2k+1}^{2l-2})(M^{-1},\sigma^-({\epsilon})) \tilde{\lambda}^k$.
 \end{enumerate}
 Here we let $\tilde{c}_{0}^{0}=1$ and $\tilde{c}_{k}^{n}=0$ for any integer $k>n$.
\end{proposition}
\begin{proof}
 The statements hold for $l=1$ by the following identities. 
 \begin{equation*}
\begin{split}
\tilde{c}_{1}^{2}&=\epsilon_1M^{-\epsilon_2}=c_{1}^{2}\\
\tilde{d}_{1}^{2}&=\epsilon_2M^{\epsilon_1}=d_{1}^{2}(M^{-1},\epsilon)\\
\tilde{c}_{0}^{2} +\tilde{c}_{2}^{2} \tilde{\lambda}&
=M^{-\epsilon_1+\epsilon_2}+\epsilon_1\epsilon_2\tilde{\lambda}
=M^{-\epsilon_1+\epsilon_2}+\epsilon_1\epsilon_2(\lambda+z^2)\\
&=M^{-\epsilon_1+\epsilon_2}+\epsilon_1\epsilon_2\lambda+z\epsilon_1z\epsilon_2\\
&=M^{-\epsilon_1+\epsilon_2}+\epsilon_1\epsilon_2\lambda+(M^{\epsilon_1}-M^{-\epsilon_1})(M^{\epsilon_2}-M^{-\epsilon_2})\\
&=M^{\epsilon_1+\epsilon_2}-M^{\epsilon_1-\epsilon_2}+M^{-\epsilon_1-\epsilon_2}+\epsilon_1\epsilon_2\lambda\\
c_{0}^{2}+ c_{2}^{2} \lambda-z c_{1}^{2}&
=M^{\epsilon_1+\epsilon_2}+\epsilon_1\epsilon_2\lambda-z \epsilon_1M^{-\epsilon_2}=M^{\epsilon_1+\epsilon_2}+\epsilon_1\epsilon_2 \lambda-(M^{\epsilon_1}-M^{-\epsilon_1})M^{-\epsilon_2}\\
&=M^{\epsilon_1+\epsilon_2}+\epsilon_1\epsilon_2 \lambda-z \epsilon_1M^{-\epsilon_2}\\
&=M^{\epsilon_1+\epsilon_2}-M^{\epsilon_1-\epsilon_2}+M^{-\epsilon_1-\epsilon_2}+\epsilon_1\epsilon_2\lambda\\
d_{0}^{2}&=M^{-\epsilon_1-\epsilon_2}=M^{-\epsilon_1-\epsilon_2}(\tilde{c}_{0}^{0}-0)\tilde{\lambda}^0\\
\end{split}
\end{equation*}
 Now suppose that the statements hold for all $l<m$ to proceed by induction on $l$. 
By Lemma \ref{lem-2} and the induction hypothesis
 \begin{equation*}
\begin{split}
\displaystyle\sum_{k=0}^{m-1} \tilde{c}_{2k+1}^{2m} \tilde{\lambda}^k
&=\displaystyle\sum_{k=0}^{m-1} 
(\tilde{c}_{2k+1}^{2m-2}M^{\epsilon_{2m-1}-\epsilon_{2m}}+\epsilon_{2m-1}\tilde{c}_{2k}^{2m-2}M^{-\epsilon_{2m}})\tilde{\lambda}^{k}\\
&=\displaystyle\sum_{k=0}^{m-1} 
\tilde{c}_{2k+1}^{2m-2}M^{\epsilon_{2m-1}-\epsilon_{2m}}\tilde{\lambda}^{k}
+\epsilon_{2m-1}M^{-\epsilon_{2m}}\displaystyle\sum_{k=0}^{m-1} \tilde{c}_{2k}^{2m-2}\tilde{\lambda}^{k}\\
&=\displaystyle\sum_{k=0}^{m-1}c_{2k+1}^{2m-2}M^{\epsilon_{2m-1}-\epsilon_{2m}}\lambda^k
+\epsilon_{2m-1}M^{-\epsilon_{2m}}\displaystyle\sum_{k=0}^{m-1}(c_{2k}^{2m-2}-zc_{2k+1}^{2m-2})\lambda^{k}
\\
&=
\displaystyle\sum_{k=0}^{m-1}(c_{2k+1}^{2m-2}(M^{\epsilon_{2m-1}-\epsilon_{2m}}-z\epsilon_{2m-1}M^{-\epsilon_{2m}})
+\epsilon_{2m-1}c_{2k}^{2m-2}M^{-\epsilon_{2m}})\lambda^{k}\\
&=
\displaystyle\sum_{k=0}^{m-1}(c_{2k+1}^{2m-2}M^{-\epsilon_{2m-1}-\epsilon_{2m}}
+\epsilon_{2m-1}c_{2k}^{2m-2}M^{-\epsilon_{2m}})\lambda^{k}\\
&=\displaystyle\sum_{k=0}^{m-1} 
c_{2k+1}^{2m} \lambda^k\\
\end{split}
\end{equation*} 
which proves (i). The proofs of (ii), (iii), and (iv) are similar.
\end{proof}

\section{Non-abelian $SL(2,\bc)$-representations of kmot groups} %=============================================
In this section, we describe the 2-variable Riley polynomials using $c_k^n, d_k^n$ defined in the previous section, and obtain several identities for the Alexander polynomials.
\begin{definition}
A $(2-bridge)$ $kmot$ $group$ is a group $G$ with a presentation  
\begin{equation}\label{kmot-group}
G=\langle x,y \,|\,wx=yw \rangle,
\end{equation}
where $w$ is of the form
\begin{equation}\label{presentation1}
w=x^{\epsilon_1}y^{\epsilon_2}x^{\epsilon_3}y^{\epsilon_4}\cdots x^{\epsilon_{\alpha-2}}y^{\epsilon_{\alpha-1}} \,\,(\alpha \,\text{is an odd integer}\,\geq 3)
\end{equation}
and  $\epsilon_i=\epsilon_{\alpha-i}=\pm 1$ for $i=1, \cdots, \alpha-1$.
\end{definition}
We write ${\bf \epsilon}=(\epsilon_1,\epsilon_2,\cdots,\epsilon_{\alpha-1})$ and $G=G( \epsilon)=G(\epsilon_1,\epsilon_2,\cdots,\epsilon_{\alpha-1})$. 
Let $F_{x,y}=\langle x,y  \rangle$ be a free group on two generators $x$ and $y$ 
and 
$\rho : F_{x,y}  \rightarrow SL(2,\bc)$ with
\begin{equation}\label{nonabel-rep}
\rho(x)=\begin{pmatrix}
M & 1\\
  0 & \frac{1}{M} 
 \end{pmatrix} \quad \text{and}\quad
\rho(y)=\begin{pmatrix}
M & 0\\
  \lambda& \frac{1}{M} 
 \end{pmatrix}
\end{equation}
be a group homomorpism and let
\begin{equation}\label{parabolic-rep2}
\begin{split}
W(\epsilon)&=\rho(w)=\rho(x)^{\epsilon_1}\rho(y)^{\epsilon_2}\rho(x)^{\epsilon_3}\rho(y)^{\epsilon_4}\cdots \rho(x)^{\epsilon_{\alpha-2}}\rho(y)^{\epsilon_{\alpha-1}}\\
&=\begin{pmatrix}
M & 1\\
  0 & \frac{1}{M} 
 \end{pmatrix}^{\epsilon_1}\begin{pmatrix}
M & 0\\
  \lambda& \frac{1}{M} 
 \end{pmatrix}^{\epsilon_2}\begin{pmatrix}
M & 1\\
  0 & \frac{1}{M} 
 \end{pmatrix}^{\epsilon_3}\cdots \begin{pmatrix}
M & 1\\
  0 & \frac{1}{M} 
 \end{pmatrix}^{\epsilon_{\alpha-2}}\begin{pmatrix}
M & 0\\
  \lambda& \frac{1}{M} 
 \end{pmatrix}^{\epsilon_{\alpha-1}}.
\end{split}
\end{equation}
Then the ($i,j$)-element $W_{ij}$ of  $W(\epsilon)=W(\epsilon_1,\epsilon_2,\cdots,\epsilon_{\alpha-1})$ has  an explicit formula for $i,j=1,2$ as follows.
\begin{proposition}\label{W-eq}
% Let $l=\frac{\alpha-1}{2}$. Then 
$$ W_{11}=\sum_{k=0}^{\frac{\alpha-1}{2}} c_{2k}^{\alpha-1} \lambda^k, \,\,W_{12}= \displaystyle\sum_{k=0}^{\frac{\alpha-3}{2}} c_{2k+1}^{\alpha-1} \lambda^k, \,\,W_{21}=\lambda W_{12}, \,\,W_{22}=\displaystyle\sum_{k=0}^{\frac{\alpha-3}{2}} d_{2k}^{\alpha-1} \lambda^k, $$ that is, 
$$W(\epsilon)=
\begin{pmatrix}
\displaystyle\sum_{k=0}^{\frac{\alpha-1}{2}} c_{2k}^{\alpha-1} \lambda^k& \displaystyle\sum_{k=0}^{\frac{\alpha-3}{2}} c_{2k+1}^{\alpha-1} \lambda^k\\
\lambda\displaystyle\sum_{k=0}^{\frac{\alpha-3}{2}} d_{2k+1}^{\alpha-1} \lambda^k&  \displaystyle\sum_{k=0}^{\frac{\alpha-3}{2}} d_{2k}^{\alpha-1} \lambda^k
 \end{pmatrix}=\begin{pmatrix}
\displaystyle\sum_{k=0}^{\frac{\alpha-1}{2}} c_{2k}^{\alpha-1} \lambda^k& \displaystyle\sum_{k=0}^{\frac{\alpha-3}{2}} c_{2k+1}^{\alpha-1} \lambda^k\\
\lambda\displaystyle\sum_{k=0}^{\frac{\alpha-3}{2}} c_{2k+1}^{\alpha-1} \lambda^k&  \displaystyle\sum_{k=0}^{\frac{\alpha-3}{2}} d_{2k}^{\alpha-1} \lambda^k
 \end{pmatrix}.$$
\end{proposition}
\begin{proof}
%Firstly, $W_{21}=\lambda W_{12}$  is easily proved  from the symmetry of $\epsilon_i$ by the  induction argument. 
%Now we prove the explicit formulas for $W_{11}, W_{12}, W_{22}$ without using the symmetry of $\epsilon_i$.
The formulas hold for $\alpha=3$ because
\begin{equation*}
\begin{split}
W(\epsilon_1, \epsilon_2)&=\rho(x)^{\epsilon_1}\rho(y)^{\epsilon_2}=\begin{pmatrix}
M & 1\\
  0 & \frac{1}{M} 
 \end{pmatrix}^{\epsilon_1}\begin{pmatrix}
M & 0\\
  \lambda& \frac{1}{M} 
 \end{pmatrix}^{\epsilon_2}=\begin{pmatrix}
M^{\epsilon_1} & \epsilon_1\\
  0 & M^{-\epsilon_1}
 \end{pmatrix}\begin{pmatrix}
M^{\epsilon_2} & 0\\
  \epsilon_2\lambda& M^{-\epsilon_2} 
 \end{pmatrix}\\
&
=\begin{pmatrix}
M^{\epsilon_1+\epsilon_2}+\lambda \epsilon_1\epsilon_2 & \epsilon_1M^{-\epsilon_2}\\
  \epsilon_2\lambda M^{-\epsilon_1}& M^{-\epsilon_1-\epsilon_2} 
 \end{pmatrix}
=\begin{pmatrix}
c_{0}^{2} +c_{2}^{2} \lambda& c_{1}^{2} \\
\lambda d_{1}^{2} &  d_{0}^{2} 
 \end{pmatrix}
=\begin{pmatrix}
c_{0}^{2} +c_{2}^{2} \lambda& c_{1}^{2} \\
\lambda c_{1}^{2} &  d_{0}^{2} 
 \end{pmatrix}.
\end{split}
\end{equation*}
Suppose that the formulas  hold for all $\alpha<2l$. Then
\begin{equation*}
\begin{split}
W(\epsilon_1, \cdots,\epsilon_{2l})&= W(\epsilon_1, \cdots,\epsilon_{2l-2}) \begin{pmatrix}
M^{\epsilon_{2l-1}+\epsilon_{2l}}+\lambda \epsilon_{2l-1}\epsilon_{2l} & \epsilon_{2l-1}M^{-\epsilon_{2l}}\\
  \epsilon_{2l}\lambda M^{-\epsilon_{2l-1}}& M^{-\epsilon_{2l-1}-\epsilon_{2l}}
 \end{pmatrix}\\
 &=\begin{pmatrix}
\displaystyle\sum_{k=0}^{l-1} c_{2k}^{2l-2} \lambda^k& \displaystyle\sum_{k=0}^{l-2} c_{2k+1}^{2l-2} \lambda^k\\
\lambda\displaystyle\sum_{k=0}^{l-2} d_{2k+1}^{2l-2} \lambda^k&  \displaystyle\sum_{k=0}^{l-2} d_{2k}^{2l-2} \lambda^k
 \end{pmatrix}
 \begin{pmatrix}
M^{\epsilon_{2l-1}+\epsilon_{2l}}+\lambda \epsilon_{2l-1}\epsilon_{2l} & \epsilon_{2l-1}M^{-\epsilon_{2l}}\\
  \epsilon_{2l}\lambda M^{-\epsilon_{2l-1}}& M^{-\epsilon_{2l-1}-\epsilon_{2l}}
 \end{pmatrix}.\\
\end{split}
\end{equation*}
We have 
\begin{equation*}
\begin{split}
W(\epsilon_1, \cdots,\epsilon_{2l})_{11}&=  \displaystyle\sum_{k=0}^{l-1} c_{2k}^{2l-2}M^{\epsilon_{2l-1}+\epsilon_{2l}} \lambda^k +\epsilon_{2l-1}\epsilon_{2l} \displaystyle\sum_{k=0}^{l-1}c_{2k}^{2l-2}\lambda^{k+1} +\epsilon_{2l}\displaystyle\sum_{k=0}^{l-2} c_{2k+1}^{2l-2}M^{-\epsilon_{2l-1}}\lambda^{k+1} \\
&=   \displaystyle\sum_{k=0}^{l-1} c_{2k}^{2l-2}M^{\epsilon_{2l-1}+\epsilon_{2l}} \lambda^k +\epsilon_{2l-1}\epsilon_{2l} \displaystyle\sum_{k=1}^{l}c_{2k-2}^{2l-2}\lambda^{k} +\epsilon_{2l}\displaystyle\sum_{k=1}^{l-1} c_{2k-1}^{2l-2}M^{-\epsilon_{2l-1}}\lambda^{k} \\
&= M^{\sum_i\epsilon_i} +\epsilon_{2l-1}\epsilon_{2l} c_{2l-2}^{2l-2} \lambda^l+  \displaystyle\sum_{k=1}^{l-1}( c_{2k}^{2l-2}M^{\epsilon_{2l-1}+\epsilon_{2l}}  +\epsilon_{2l-1}\epsilon_{2l} c_{2k-2}^{2l-2} +\epsilon_{2l} c_{2k-1}^{2l-2}M^{-\epsilon_{2l-1}})\lambda^{k}\\
&= c_{0}^{2l}+ c_{2l}^{2l} \lambda^l+  \displaystyle\sum_{k=1}^{l-1}( c_{2k}^{2l-2}M^{\epsilon_{2l-1}+\epsilon_{2l}}  +\epsilon_{2l-1}\epsilon_{2l} c_{2k-2}^{2l-2} +\epsilon_{2l} c_{2k-1}^{2l-2}M^{-\epsilon_{2l-1}})\lambda^{k},\\
%&= M^{\sum_i\epsilon_i}+  \displaystyle\sum_{k=1}^{l-1}c_{2k}^{2l}\lambda^{k} +c_{2l}^{2l} \lambda^l\\
\end{split}
\end{equation*}
\begin{equation*}
\begin{split}
W(\epsilon_1, \cdots,\epsilon_{2l})_{12}&=   \displaystyle\sum_{k=0}^{l-1} \epsilon_{2l-1}c_{2k}^{2l-2}M^{-\epsilon_{2l}}\lambda^k +  \displaystyle\sum_{k=0}^{l-2} c_{2k+1}^{2l-2}M^{-\epsilon_{2l-1}-\epsilon_{2l}}\lambda^k\\
&=\epsilon_{2l-1}c_{2l-2}^{2l-2}M^{-\epsilon_{2l}}\lambda^{l-1}+   \displaystyle\sum_{k=0}^{l-2}( c_{2k+1}^{2l-2}M^{-\epsilon_{2l-1}-\epsilon_{2l}}+\epsilon_{2l-1}c_{2k}^{2l-2}M^{-\epsilon_{2l}})\lambda^k\\
&= c_{2l-1}^{2l}\lambda^{l-1}+\displaystyle\sum_{k=0}^{l-2}( c_{2k+1}^{2l-2}M^{-\epsilon_{2l-1}-\epsilon_{2l}}+\epsilon_{2l-1}c_{2k}^{2l-2}M^{-\epsilon_{2l}})\lambda^k,\\
\end{split}
\end{equation*}
\begin{equation*}
\begin{split}
W(\epsilon_1, \cdots,\epsilon_{2l})_{21}&=   \lambda\left(\displaystyle\sum_{k=0}^{l-2} d_{2k+1}^{2l-2}M^{\epsilon_{2l-1}+\epsilon_{2l}}\lambda^k+ \epsilon_{2l-1}\epsilon_{2l}d_{2k+1}^{2l-2}\lambda^{k+1}\right)
+\displaystyle\sum_{k=0}^{l-2}\epsilon_{2l}d_{2k}^{2l-2}M^{-\epsilon_{2l-1}}\lambda^{k+1}
\\
&=   \lambda\left( \epsilon_{2l-1}\epsilon_{2l}d_{2l-3}^{2l-2}\lambda^{l-1}+\displaystyle\sum_{k=0}^{l-2} (d_{2k+1}^{2l-2}M^{\epsilon_{2l-1}+\epsilon_{2l}}+\epsilon_{2l}d_{2k}^{2l-2}M^{-\epsilon_{2l-1}}+ \epsilon_{2l-1}\epsilon_{2l}d_{2k-1}^{2l-2})\lambda^{k}\right)\\
&=   \lambda\left(d_{2l-1}^{2l}\lambda^{l-1}+\displaystyle\sum_{k=0}^{l-2} (d_{2k+1}^{2l-2}M^{\epsilon_{2l-1}+\epsilon_{2l}}+\epsilon_{2l}d_{2k}^{2l-2}M^{-\epsilon_{2l-1}}+ \epsilon_{2l-1}\epsilon_{2l}d_{2k-1}^{2l-2})\lambda^{k}\right)\\
\end{split}
\end{equation*}
and
\begin{equation*}
\begin{split}
W(\epsilon_1, \cdots,\epsilon_{2l})_{22}&= \displaystyle\sum_{k=0}^{l-2} ( d_{2k}^{2l-2}M^{-\epsilon_{2l-1}-\epsilon_{2l}}\lambda^k+\epsilon_{2l-1}d_{2k+1}^{2l-2}M^{-\epsilon_{2l}}\lambda^{k+1})\\
&= \epsilon_{2l-1}d_{2l-3}^{2l-2}M^{-\epsilon_{2l}}\lambda^{l-1}+\displaystyle\sum_{k=0}^{l-2}  (d_{2k}^{2l-2}M^{-\epsilon_{2l-1}-\epsilon_{2l}}+\epsilon_{2l-1}d_{2k-1}^{2l-2}M^{-\epsilon_{2l}})\lambda^{k}\\
&= d_{2l}^{2l-1}\lambda^{l-1}+ \displaystyle\sum_{k=0}^{l-2}  (d_{2k}^{2l-2}M^{-\epsilon_{2l-1}-\epsilon_{2l}}+\epsilon_{2l-1}d_{2k-1}^{2l-2}M^{-\epsilon_{2l}})\lambda^{k}.   \\
\end{split}
\end{equation*}
Hence by Lemma \ref{symmetric} and \ref{lem-2}  we have 
 $$W(\epsilon_1, \cdots,\epsilon_{2l})_{11}=c^{2l}_0+c_{2l}^{2l} \lambda^l+  \displaystyle\sum_{k=1}^{l-1}c_{2k}^{2l}\lambda^{k}=\displaystyle\sum_{k=0}^{l}c_{2k}^{2l}\lambda^{k},$$
$$W(\epsilon_1, \cdots,\epsilon_{2l})_{12}=c_{2l-1}^{2l}\lambda^{l-1}+\displaystyle\sum_{k=0}^{l-2}c_{2k+1}^{2l} \lambda^k=\displaystyle\sum_{k=0}^{l-1} c_{2k+1}^{2l} \lambda^k,$$
$$W(\epsilon_1, \cdots,\epsilon_{2l})_{21}=\lambda(d_{2l-1}^{2l}\lambda^{l-1}+\displaystyle\sum_{k=0}^{l-2}d_{2k+1}^{2l} \lambda^k)=\lambda\displaystyle\sum_{k=0}^{l-1} d_{2k+1}^{2l} \lambda^k=\lambda\displaystyle\sum_{k=0}^{l-1} c_{2k+1}^{2l} \lambda^k,$$
and
$$W(\epsilon_1, \cdots,\epsilon_{2l})_{22}=d_{2l}^{2l-1}\lambda^{l-1}+\displaystyle\sum_{k=0}^{l-2}d_{2k}^{2l} \lambda^k
=\displaystyle\sum_{k=0}^{l-1}d_{2k}^{2l} \lambda^k, $$
which proves that  the formulas hold for $\alpha=2l+1$ also.
\end{proof}

Using the two identities $W_{21}=\lambda W_{12}$ and
$\begin{pmatrix}
W_{11} &W_{12}\\
W_{21}&  W_{22}
 \end{pmatrix}
 \begin{pmatrix}
M & 1\\
  0 & \frac{1}{M} 
 \end{pmatrix} =\begin{pmatrix}
M & 0\\
  \lambda& \frac{1}{M} 
 \end{pmatrix}
 \begin{pmatrix}
W_{11} &W_{12}\\
W_{21}&  W_{22}
 \end{pmatrix}$, the following Riley's result is not difficult to prove, but here we use the quandle argument to provide another proof.
\begin{proposition} [\cite{Riley2}]\label{Riley poly}
$\rho$ defines a non-abelian representation on the kmot group $G(\epsilon)$ if and only if  $W_{11}-zW_{12}=0$.
\end{proposition}
\begin{proof}
Let $W=W(\epsilon)$.
The relation $W\rho(x)W^{-1}=\rho(w)\rho(x)\rho(w)^{-1}=\rho(y)$ with (\ref{nonabel-rep}) is equivalent to  $\tilde{W}\sim \tilde{b}$ by (\ref{b}), where 
$ b=\begin{pmatrix}
z & 1\\
\lambda& \frac{\lambda+1}{z}
 \end{pmatrix},$ 
and thus 
\begin{equation}\label{w-b}
\begin{pmatrix}
W_{11} & W_{11}-zW_{12}\\
W_{21}&  W_{21}-zW_{22}
 \end{pmatrix}=\tilde{W}\sim \tilde{b}= \begin{pmatrix}
z & 0\\
\lambda& -1
 \end{pmatrix}.
\end{equation}
This equation  is equivalent to 
$$
\begin{pmatrix}
W_{11} & W_{11}-zW_{12}\\
\lambda W_{12}&  W_{21}-zW_{22}
 \end{pmatrix}=\begin{pmatrix}
z & 0\\
\lambda& -1
 \end{pmatrix}
\begin{pmatrix}
t & 0\\
0& t^{-1}
 \end{pmatrix}
=\begin{pmatrix}
zt & 0\\
\lambda t& - t^{-1}
 \end{pmatrix} \quad \text{for some}\quad  t\in \mathbb C,$$
which holds when $W_{11}-zW_{12}=0$ and $z+W_{11}( W_{21}-zW_{22})=0$. But 
$$W_{11}-zW_{12} \,|\,z+W_{11}( W_{21}-zW_{22})$$ since if $W_{11}-zW_{12}=0$ then 
\begin{equation*}
\begin{split}
z+W_{11}(W_{21}-zW_{22})&= z+zW_{12}(W_{21}-zW_{22})=z+z(W_{11}W_{22}-1)-z^2W_{12}W_{22}\\
&=zW_{22}(W_{11}-zW_{12})=0. 
\end{split}
\end{equation*}
This completes the proof.
\end{proof}
\begin{definition}
We call the 2-variable polynomial  $W_{11}-zW_{12}$ in $\mathbb Z[M^{\pm 1},\lambda]$ the $Riley$ $polynomial$ of the kmot group of $G(\epsilon)$. 
\end{definition}
These Riley polynomials are related to the normalized Alexander polynomials as follows:
\begin{corollary}\label{Alexander1}
The Alexander polynomial $\bigtriangleup(t)$ of a kmot group $G( \epsilon)$ satisfies the following.
$$\bigtriangleup(M^2)=\bigtriangleup(M^{-2})=M^{\sigma}-z
\sum_{i_1:\text {odd}} \epsilon_{i_1} M^{\widehat{i_1}},\quad \sigma=\displaystyle \sum_{i=1}^{\alpha-1} \epsilon_i $$
\end{corollary}
\begin{proof}
By Riley \cite{Riley2}, the Alexander polynomial is obtained when $\lambda=0$ and 
$\bigtriangleup(M^2)$ is  $W_{11}(M,0)-z W_{12}(M,0)=W_{11}(M^{-1},0)-z W_{12}(M^{-1},0)$ up to $\mathbb Z[M^2, M^{-2}]$-unit multiple. Hence 
we have $$\bigtriangleup(M^2)=W_{11}(M,0)-z W_{12}(M,0)=M^{\sigma}-z
\sum_{i_1:\text {odd}}^{\wedge} \epsilon_{i_1} M^{\widehat{i_1}}=M^{\sigma}-z
\sum_{i_1:\text {odd}} \epsilon_{i_1} M^{\widehat{i_1}}.$$
Note that $\displaystyle\sum_{i_1:\text {odd}}^{\wedge}=\sum_{i_1:\text {odd}}$. 
\end{proof}

From Corollary \ref{Alexander1}, we have the following formulas.
\begin{corollary}\label{Min}
Let  $\bigtriangleup(t)$ be the Alexander polynomial of a kmot group $G( \epsilon)$. % and $\nu_k=1+\sum_{i=1}^{\alpha-1}\epsilon_{k+i}$.
 Then we have the following.
\begin{enumerate}
\item [\rm (i)]
 $ \bigtriangleup(t)=t^{\frac{\sigma}{2}}+ (t^{-\frac{1}{2}}-t^{\frac{1}{2}})\sum_{k:odd}\epsilon_{k}t^{\frac{1}{2}(\epsilon_{1}+\cdots+\epsilon_{k-1}-(\epsilon_{k+1}+\cdots+\epsilon_{\alpha-1}))}=\bigtriangleup(t^{-1})$
\item [\rm (ii)] $t^{\frac{\sigma}{2}}\bigtriangleup(t)=1-t^{\epsilon_{1}}+t^{\epsilon_{1}+\epsilon_{2}}-t^{\epsilon_{1}+\epsilon_{2}+\epsilon_{3}}+\cdots+
(-1)^{k}t^{\epsilon_{1}+\cdots+\epsilon_{k}}+\cdots
t^{\epsilon_{1}+\cdots+\epsilon_{\alpha-1}} $
\end{enumerate}
\end{corollary}
\begin{proof} 
\begin{enumerate}
\item [\rm (i)]
The identity comes straight from:
\begin{equation*}
\begin{split}
\bigtriangleup(M^2)&=M^{\sigma}-z\sum_{i_1:\text {odd}} \epsilon_{i_1} M^{\widehat{i_1}}=M^{\sigma}-(M-\frac{1}{M})\sum_{i_1:\text {odd}} \epsilon_{i_1} M^{\widehat{i_1}}
\end{split}
\end{equation*}
\item [\rm (ii)]
Since $\widehat{k}=-\sigma+\epsilon_{k}+2\sum_{i=1}^{k-1}\epsilon_{i}$ (recall Definition \ref{cd} for $\widehat{k}$), 
\begin{equation*}
\begin{split}
\bigtriangleup(M^2)&=M^{\sigma}-z\sum_{i_1:\text {odd}} \epsilon_{i_1} M^{\widehat{i_1}}\\
&=M^{\sigma}-M^{-\sigma}(z\epsilon_{1}M^{\epsilon_{1}}+
z\epsilon_{3}M^{2(\epsilon_{1}+\epsilon_{2})+\epsilon_{3}}
+\cdots+z\epsilon_{\alpha-2}M^{2(\epsilon_{1}+\cdots+\epsilon_{\alpha-3})+\epsilon_{\alpha-2}})\\
&=M^{\sigma}-M^{-\sigma}((M^{\epsilon_{1}}-M^{-\epsilon_{1}})M^{\epsilon_{1}}
+\cdots+(M^{\epsilon_{\alpha-2}}-M^{-\epsilon_{\alpha-2}})M^{2(\epsilon_{1}+\cdots+\epsilon_{\alpha-3})+\epsilon_{\alpha-2}})\\
&=M^{\sigma}-M^{-\sigma}(-1+M^{2\epsilon_{1}}-M^{2(\epsilon_{1}+\epsilon_{2})}+M^{2(\epsilon_{1}+\epsilon_{2}+\epsilon_{3})}\cdots+M^{2(\epsilon_{1}+\cdots+\epsilon_{\alpha-2})})\\
&=-M^{-\sigma}(-1+M^{2\epsilon_{1}}-M^{2(\epsilon_{1}+\epsilon_{2})}+M^{2(\epsilon_{1}+\epsilon_{2}+\epsilon_{3})}\cdots-M^{2(\epsilon_{1}+\cdots+\epsilon_{\alpha-1})}),
\end{split}
\end{equation*}
and the identity is obtained by multiplying $M^{\sigma}=t^{\frac{\sigma}{2}}$.
\end{enumerate}
%we have the following formula of Minkus in \cite{Minkus} (see also \cite{Hoste} or \cite{HS2}) by multiplying $M^{\sigma}=t^{\frac{\sigma}{2}}$ if we let $t=M^2$.
%$$\bigtriangleup(t)=1-t^{\epsilon_{1}}+t^{\epsilon_{1}+\epsilon_{2}}-t^{\epsilon_{1}+\epsilon_{2}+\epsilon_{3}}+\cdots+
%(-1)^{k+1}t^{\epsilon_{1}+\cdots+\epsilon_{k}}+\cdots
%t^{\epsilon_{1}+\cdots+\epsilon_{\alpha-1}} $$
%We can also show that our formula  is equal to 
%Fukuhara's formula in \cite{Fukuhara},
%$$ \bigtriangleup(t)=\frac{1}{2}(t^{-\frac{\sigma}{2}}+t^{\frac{\sigma}{2}})- \frac{1}{4}(t^{-\frac{1}{2}}-t^{\frac{1}{2}})\sum_{k=1}^{\alpha-1}(-1)^k\epsilon_{k}(t^{-\frac{\nu_k}{2}}-t^{\frac{\nu_k}{2}}) $$
%where $\nu_k=1+\sum_{i=1}^{\alpha-1}\epsilon_{k+i}$.
\end{proof}
Note that the right side of (ii) is exactly Minkus' formula in \cite{Minkus} (see also \cite{Hoste} or \cite{HS2}), and (i) is related to Fukuhara's formula in \cite{Fukuhara} as follows: 
If we let 
\begin{equation}\label{Fuk}
  \epsilon_0=1, \epsilon_{\alpha+i}= -\epsilon_i,
\end{equation}
 $-\widehat{k}$ is equal to $\nu_k=1+\sum_{i=1}^{\alpha-1}\epsilon_{k+i}$, and (i) is equivalent to 
 \begin{equation}\label{AP}
  \bigtriangleup(t)=t^{\frac{\sigma}{2}}+ (t^{-\frac{1}{2}}-t^{\frac{1}{2}})\sum_{\shortstack{$\scriptstyle k:\text{odd} $\\$\scriptstyle  1\leq k\leq \alpha-1$}}\epsilon_{k}t^{-\frac{\nu_k}{2}}=\bigtriangleup(t^{-1}).
  \end{equation}
Fukuhara defined  $\nu_k$ in \cite{Fukuhara} for the sequence $\epsilon_i=(-1)^{[i\frac{\beta}{\alpha}]}$ corresponding to a 2-bridge knot $K=S(\alpha,\beta)$ to get a formula for the  Alexander polynomial of $K$, and in this case (\ref{Fuk}) is automatically satisfied. (Note that the knot group $G(K)$ is a kmot group.) 
\begin{corollary}[\cite{Fukuhara}]
Let  $\bigtriangleup(t)$ be the Alexander polynomial of a kmot group $G(\epsilon_1,\epsilon_2,\cdots,\epsilon_{\alpha-1})$. 
Let 
$\epsilon_0=1, \epsilon_{\alpha+i}= -\epsilon_i$ and $\nu_k=1+\sum_{i=1}^{\alpha-1}\epsilon_{k+i}$ for integers $i,k\in\{1, \cdots, \alpha-1\}$. Then we have
$$ \bigtriangleup(t)=\frac{1}{2}(t^{-\frac{\sigma}{2}}+t^{\frac{\sigma}{2}})- \frac{1}{4}(t^{-\frac{1}{2}}-t^{\frac{1}{2}})\sum_{k=1}^{\alpha-1}(-1)^k\epsilon_{k}(t^{-\frac{\nu_k}{2}}-t^{\frac{\nu_k}{2}}) $$
\end{corollary}
\begin{proof}
Since $\epsilon_{\alpha-k}= \epsilon_k$ and $\nu_k=-\nu_{\alpha-k}$, 
\begin{equation*}
\begin{split}
\bigtriangleup(t)&=t^{\frac{\sigma}{2}}+(t^{-\frac{1}{2}}-t^{\frac{1}{2}})\sum_{k:odd}\epsilon_{k}t^{-\frac{\nu_k}{2}}
%=t^{\frac{\sigma}{2}}+(t^{-\frac{1}{2}}-t^{\frac{1}{2}})\sum_{k:odd}\epsilon_{\alpha-k}t^{\frac{\nu_{\alpha-k}}{2}}\\
  =t^{\frac{\sigma}{2}}+(t^{-\frac{1}{2}}-t^{\frac{1}{2}})\sum_{k:even}\epsilon_{k}t^{\frac{\nu_{k}}{2}},
\end{split}
\end{equation*}
\begin{equation*}
\begin{split}
\bigtriangleup(t^{-1})&=t^{-\frac{\sigma}{2}}-(t^{-\frac{1}{2}}-t^{\frac{1}{2}})\sum_{k:odd}\epsilon_{k}t^{\frac{\nu_k}{2}}=t^{-\frac{\sigma}{2}}-(t^{-\frac{1}{2}}-t^{\frac{1}{2}})\sum_{k:even}\epsilon_{k}t^{-\frac{\nu_{k}}{2}}.
\end{split}
\end{equation*}
Hence the identity is proved as follows:
\begin{equation*}
\begin{split}
4\bigtriangleup(t)&=2(\bigtriangleup(t)+\bigtriangleup(t^{-1}))\\
&=2(t^{-\frac{\sigma}{2}}+t^{\frac{\sigma}{2}})-(t^{-\frac{1}{2}}-t^{\frac{1}{2}})\sum_{k=1}^{\alpha-1}(-1)^k\epsilon_{k}(t^{-\frac{\nu_k}{2}}-t^{\frac{\nu_k}{2}})
\end{split}
\end{equation*}
\end{proof}

Now let's look at the definition of $W^*$, which we will use later.
%Now we turn to the definition of $W^*$, which will be used later.
Let $w^*$ be the element  obtained from $w$ by exchanging $x$ and $y$, that is, 
\begin{equation}\label{presentation2}
w^*=y^{\epsilon_1}x^{\epsilon_2}y^{\epsilon_3}x^{\epsilon_4}\cdots y^{\epsilon_{\alpha-2}}x^{\epsilon_{\alpha-1}},
\end{equation}
 and  $\bar{W}_{ij}(M,\lambda):=W_{ij}(M^{-1},\lambda)$.
Then we have
\begin{lemma}[\cite{HS}]\label{w-star}
 $ W^*=\rho(w^*)=\begin{pmatrix}
\bar{W}_{22} & \bar{W}_{12} \\
\lambda \bar{W}_{12}& \bar{W}_{11}
 \end{pmatrix}$.
\end{lemma}
\begin{proof}
The identity holds for $\alpha=3$ because
\begin{equation*}
\begin{split}
W(\epsilon_1, \epsilon_2)&=\rho(x)^{\epsilon_1}\rho(y)^{\epsilon_2}=\begin{pmatrix}
M^{\epsilon_1} & \epsilon_1\\
  0 & M^{-\epsilon_1}
 \end{pmatrix}\begin{pmatrix}
M^{\epsilon_2} & 0\\
  \epsilon_2\lambda& M^{-\epsilon_2} 
 \end{pmatrix}=\begin{pmatrix}
M^{\epsilon_1+\epsilon_2}+\lambda \epsilon_1\epsilon_2 & \epsilon_1M^{-\epsilon_2}\\
  \epsilon_2\lambda M^{-\epsilon_1}& M^{-\epsilon_1-\epsilon_2} 
 \end{pmatrix}.
\end{split}
\end{equation*}
and
\begin{equation*}
\begin{split}
W^*(\epsilon_1, \epsilon_2)&=\rho(y)^{\epsilon_1}\rho(x)^{\epsilon_2}=\begin{pmatrix}
M^{\epsilon_2} & 0\\
  \epsilon_2\lambda& M^{-\epsilon_2} 
 \end{pmatrix}\begin{pmatrix}
M^{\epsilon_1} & \epsilon_1\\
  0 & M^{-\epsilon_1}
 \end{pmatrix}=\begin{pmatrix}
M^{\epsilon_1+\epsilon_2}& \epsilon_1M^{\epsilon_2}\\
  \epsilon_2\lambda M^{\epsilon_1}& M^{-\epsilon_1-\epsilon_2} +\lambda \epsilon_1\epsilon_2
 \end{pmatrix}.
\end{split}
\end{equation*}
Assume that the identity  holds for all $\alpha<2n$. Then
\begin{equation*}
\begin{split}
W^*(\epsilon_1, \cdots,\epsilon_{2n})&= \begin{pmatrix}
M^{\epsilon_1} & 0\\
  \epsilon_1\lambda & M^{-\epsilon_1}
 \end{pmatrix}W(\epsilon_2, \cdots,\epsilon_{2n-1}) \begin{pmatrix}
M^{\epsilon_{2n}} & \epsilon_{2n}\\
  0& M^{-\epsilon_{2n}} 
 \end{pmatrix}\\
W(\epsilon_1, \cdots,\epsilon_{2n})&= \begin{pmatrix}
M^{\epsilon_1} & \epsilon_1\\
  0 & M^{-\epsilon_1}
 \end{pmatrix}W^*(\epsilon_2, \cdots,\epsilon_{2n-1}) \begin{pmatrix}
M^{\epsilon_{2n}} & 0\\
 \epsilon_{2n}\lambda& M^{-\epsilon_{2n}} 
 \end{pmatrix}.
\end{split}
\end{equation*}
If we let $W'=W(\epsilon_2, \cdots,\epsilon_{2n-1})$, then 
\begin{equation*}
\begin{split}
W^*(\epsilon_1, \cdots,\epsilon_{2n})&= \begin{pmatrix}
M^{\epsilon_1} & 0\\
   \epsilon_1\lambda & M^{-\epsilon_1}
 \end{pmatrix}
\begin{pmatrix}
W'_{11} & W'_{12}\\
\lambda W'_{12}&  W'_{22}
 \end{pmatrix}
\begin{pmatrix}
M^{\epsilon_{2n}} & \epsilon_{2n}\\
  0& M^{-\epsilon_{2n}} 
 \end{pmatrix}\\
&=\begin{pmatrix}
M^{\epsilon_1+\epsilon_{2n}}W'_{11}& W'_{12}M^{\epsilon_1-\epsilon_{2n}}+\epsilon_1 M^{\epsilon_{2n}} W'_{11} \\
\lambda (W'_{12}+\epsilon_1 M^{\epsilon_{2n}} W'_{11})& *
 \end{pmatrix}\\
&=\begin{pmatrix}
M^{2\epsilon_1}W'_{11}& W'_{12}+\epsilon_1 M^{\epsilon_{1}} W'_{11} \\
\lambda (W'_{12}+\epsilon_1 M^{\epsilon_{1}} W'_{11})& *
 \end{pmatrix}
\end{split}
\end{equation*}
and by the induction hypothesis 
\begin{equation*}
\begin{split}
W(\epsilon_1, \cdots,\epsilon_{2n})&= \begin{pmatrix}
M^{\epsilon_1} & \epsilon_1\\
   & M^{-\epsilon_1}
 \end{pmatrix}\begin{pmatrix}
\bar{W'}_{22} & \bar{W'}_{12} \\
\lambda \bar{W'}_{12}& \bar{W'}_{11}
 \end{pmatrix}\begin{pmatrix}
M^{\epsilon_{2n}} & 0\\
  \epsilon_{2n}\lambda& M^{-\epsilon_{2n}} 
 \end{pmatrix}\\
&=\begin{pmatrix}
* & \bar{W'}_{12}+\epsilon_1 M^{-\epsilon_1} \bar{W'}_{11} \\
\lambda (\bar{W'}_{12}+\epsilon_1 M^{-\epsilon_1} \bar{W'}_{11})& M^{-2\epsilon_1} \bar{W'}_{11}
 \end{pmatrix}.\\
\end{split}
\end{equation*}
This implies that  the identity  holds for $\alpha=2n+1$ also.

\end{proof}

\section{Symplectic quandle method}\label{QM}%=======================================================================
In this section, we use the generalized symplectic quandle structure corresponding to  ($\mathcal{D}_M$, conjugation) to study 
the 2-variable Riley polynomial of a kmot group $G(\epsilon_1,\epsilon_2,\cdots,\epsilon_{\alpha-1})$. The results of this alternative method will derive some recursive formulas for the Riley polynomials and the Alexander polynomials
%, and A-polynomials 
in the following sections.

Define $x_{i}$ for each $i=0,1,2,\cdots,\alpha-1$ as follows.
\begin{equation*}
\begin{split}
x_0&=x\\
x_1&=y^{\epsilon_{\alpha-1}}xy^{-\epsilon_{\alpha-1}}=x\rhd^{-\epsilon_{\alpha-1}}y\\
x_2&=(x^{\epsilon_{\alpha-2}}y^{\epsilon_{\alpha-1}})x(y^{-\epsilon_{\alpha-1}}x^{-\epsilon_{\alpha-2}})=x^{\epsilon_{\alpha-2}}x_1x^{-\epsilon_{\alpha-2}}=x_1\rhd^{-\epsilon_{\alpha-2}}x\\
&\cdots\\
x_{2k}&=(x^{\epsilon_{\alpha-2k}}y^{\epsilon_{\alpha-2k+1}}\cdots x^{\epsilon_{\alpha-2}}y^{\epsilon_{\alpha-1}})x(y^{-\epsilon_{\alpha-1}}x^{-\epsilon_{\alpha-2}}\cdots y^{-\epsilon_{\alpha-2k+1}}x^{-\epsilon_{\alpha-2k}})=x_{2k-1}\rhd^{-\epsilon_{\alpha-2k}}x\\
x_{2k+1}&=(y^{\epsilon_{\alpha-(2k+1)}}x^{\epsilon_{\alpha-2k}}\cdots x^{\epsilon_{\alpha-2}}y^{\epsilon_{\alpha-1}})x(y^{-\epsilon_{\alpha-1}}x^{-\epsilon_{\alpha-2}}\cdots x^{-\epsilon_{\alpha-2k}}y^{-\epsilon_{\alpha-(2k+1)}})\\
&=y^{\epsilon_{\alpha-(2k+1)}}x_{2k}y^{-\epsilon_{\alpha-(2k+1)}}=x_{2k}\rhd^{-\epsilon_{\alpha-(2k+1)}}y\\
x_{2k+2}&=x^{\epsilon_{\alpha-(2k+2)}}x_{2k+1}x^{-\epsilon_{\alpha-(2k+2)}}=x_{2k+1}\rhd^{-\epsilon_{\alpha-(2k+2)}}x\\
&\cdots\\
x_{\alpha-1}&=x^{\epsilon_{1}}\cdots y^{\epsilon_{\alpha-1}}xy^{-\epsilon_{\alpha-1}}\cdots x^{-\epsilon_{1}}=(\cdots((x\rhd^{-\epsilon_{\alpha-1}}y)\rhd^{-\epsilon_{\alpha-2}}x)\rhd^{-\epsilon_{\alpha-3}}y\cdots)\rhd^{-\epsilon_1}x
\end{split}
\end{equation*}
Note that $wxw^{-1}=y$ holds if and only if $x_{\alpha-1}=y$.
\begin{proposition}\label{kmot-rep}
Let $F_{x,y}=\langle x,y  \rangle$ be a free group on two letters $x$ and $y$, and  $\rho : F_{x,y}  \rightarrow SL(2,\bc)$ be a group homomorphism  satisfying 
$$\rho(x)=ama^{-1},\,\, \rho(y)=bmb^{-1},\,\, \rho(x)\rho(y)\neq \rho(y)\rho(x).$$
Suppose that $(\epsilon_1,\epsilon_2,\cdots,\epsilon_{\alpha-1})$ is a symmetric $\epsilon_i$-sequence, that is, $\epsilon_i=\epsilon_{\alpha-i}$ for $i=1, \cdots, \alpha-1$, and  $\lambda=2-tr(\rho(x)\rho(y)^{-1})$. 
Then there exist $F_i, G_i \in gl(2,\mathbb Z[M^{\pm 1},\tilde{\lambda}])$ such that
 $$\rho(x_i)=c_i mc_i^{-1},\quad 
  \tilde{c_i}\sim \tilde{a}F_i+\tilde{b}G_i$$ for each $i=0,1,\cdots,\alpha-1$, which satisfy the followings.
\begin{enumerate}
\item [\rm (i)] $F_0=Id, G_0=0,\,\, F_1=e^{-\epsilon_1}, G_1=-\epsilon_1[\tilde{a}, \tilde{b}]$
\item [\rm (ii)] $F_{2n}=F_{2n-1}e^{\epsilon_{2n}}-\epsilon_{2n}[\tilde{b}, \tilde{a}]G_{2n-1},\,\, G_{2n}=G_{2n-1}e^{-\epsilon_{2n}}$
\item [\rm (iii)] $F_{2n+1}=F_{2n}e^{-\epsilon_{2n+1}},  G_{2n+1}=G_{2n}e^{\epsilon_{2n+1}}-\epsilon_{2n+1}[\tilde{a}, \tilde{b}]F_{2n}$
 \end{enumerate}
\end{proposition}
\begin{proof}
%We play the quandle method developped in Proposition \ref{quandle-method} to prove all of the statements (i) to (v).
Obviously $F_0=Id, G_0=0$  because $\tilde{c_0}\sim \tilde{a}$.
From $x_1=y^{\epsilon_{\alpha-1}}xy^{-\epsilon_{\alpha-1}}=y^{\epsilon_{1}}xy^{-\epsilon_1}=x\rhd^{-\epsilon_{1}}y$ and $\rho(x_1)=c_1mc_1^{-1} $, we have
$$\tilde{c_1}\sim \tilde{a}e^{-\epsilon_1}-\epsilon_1 \tilde{b}[\tilde{a}, \tilde{b}]=\tilde{a}e^{-\epsilon_1}+ \tilde{b}(-\epsilon_1[\tilde{a}, \tilde{b}])$$
by Proposition \ref{quandle-method}. Hence (i) is proved. (ii) and (iii) also  hold for $n=0,1$, because 
\begin{equation*}
\begin{split}
\tilde{c_2}&\sim \tilde{c_1}e^{-\epsilon_2}-\epsilon_2 \tilde{a}[\tilde{c_1}, \tilde{a}]\\
&=(\tilde{a}F_1+ \tilde{b}G_1)e^{-\epsilon_2}-\epsilon_2 \tilde{a}[\tilde{a}F_1+ \tilde{b}G_1, \tilde{a}]\\
&=\tilde{a}(F_1(e^{-\epsilon_{2}}-\epsilon_{2} [\tilde{a},\tilde{a}])-\epsilon_2G_1[\tilde{b}, \tilde{a}])+ \tilde{b}(G_1e^{-\epsilon_2})\\
&=\tilde{a}(F_1e^{\epsilon_2}-\epsilon_2G_1[\tilde{b}, \tilde{a}])+ \tilde{b}(G_1e^{-\epsilon_2}),\\
\end{split}
\end{equation*}
where the last equality comes from the identity 
$$[\tilde{a},\tilde{a}]= \begin{pmatrix}
-z & 0\\
 0& z
 \end{pmatrix}=[\tilde{b},\tilde{b}], $$ 
 and the identity
\begin{equation}\label{EQ-4.1}
  e^{\mp 1}\mp
\begin{pmatrix}
-z & 0\\
 0& z
 \end{pmatrix}=e^{\pm 1}. \,\,\text{(Recall that}\,\, e=\begin{pmatrix}
M & 1\\
  0 & \frac{1}{M} 
 \end{pmatrix}.)
 \end{equation}

Now suppose that (ii) and (iii) hold for all $n<k$. 
 Then 
from $$x_{2k}=x_{2k-1}\rhd^{-\epsilon_{\alpha-2k}}x=x_{2k-1}\rhd^{-\epsilon_{2k}}x$$ and $\rho(x_{2k})=c_{2k}mc_{2k}^{-1} $, we have
\begin{equation*}
\begin{split}
\tilde{c}_{2k}&\sim \tilde{c}_{2k-1}e^{-\epsilon_{2k}}-\epsilon_{2k} \tilde{a}[\tilde{c}_{2k-1}, \tilde{a}]\\
&= (\tilde{a}F_{2k-1}+\tilde{b}G_{2k-1})e^{-\epsilon_{2k}}-\epsilon_{2k} \tilde{a}[\tilde{a}F_{2k-1}+\tilde{b}G_{2k-1}, \tilde{a}]\\
&=\tilde{a}(F_{2k-1}e^{-\epsilon_{2k}}-\epsilon_{2k} [\tilde{a},\tilde{a}]F_{2k-1}-\epsilon_{2k} [\tilde{b},\tilde{a}]G_{2k-1})+\tilde{b}G_{2k-1}e^{-\epsilon_{2k}}\\
&=\tilde{a}(F_{2k-1}(e^{-\epsilon_{2k}}-\epsilon_{2k} [\tilde{a},\tilde{a}])-\epsilon_{2k} [\tilde{b},\tilde{a}]G_{2k-1})+\tilde{b}G_{2k-1}e^{-\epsilon_{2k}}
\end{split}
\end{equation*}
Since
  $e^{-\epsilon_{2k}}-\epsilon_{2k} 
[\tilde{a},\tilde{a}]=e^{-\epsilon_{2k}}-\epsilon_{2k} 
\begin{pmatrix}
-z & 0\\
 0& z
 \end{pmatrix}=e^{\epsilon_{2k}}$ by (\ref{EQ-4.1}), we have 
$$\tilde{c}_{2k}\sim\tilde{a}(F_{2k-1}e^{\epsilon_{2k}}-\epsilon_{2k} [\tilde{b},\tilde{a}]G_{2k-1})+\tilde{b}G_{2k-1}e^{-\epsilon_{2k}},$$
which implies that
 (ii)  holds for $n=k$ also, that is, 
$$F_{2k}=F_{2k-1}e^{\epsilon_{2k}}-\epsilon_{2k} [\tilde{b},\tilde{a}]G_{2k-1},\quad G_{2k}=G_{2k-1}e^{-\epsilon_{2k}}.$$
From 
$$x_{2k+1}=x_{2k}\rhd^{-\epsilon_{\alpha-(2k+1)}}y=x_{2k}\rhd^{-\epsilon_{2k+1}}y$$ and $\rho(x_{2k+1})=c_{2k+1}mc_{2k+1}^{-1} $, we have
\begin{equation*}
\begin{split}
\tilde{c}_{2k+1}&\sim \tilde{c}_{2k}e^{-\epsilon_{2k+1}}-\epsilon_{2k+1} \tilde{b}[\tilde{c}_{2k}, \tilde{b}]\\
&= (\tilde{a}F_{2k}+\tilde{b}G_{2k})e^{-\epsilon_{2k+1}}-\epsilon_{2k+1} \tilde{b}[\tilde{a}F_{2k}+\tilde{b}G_{2k}, \tilde{b}]\\
&=\tilde{a}(F_{2k}e^{-\epsilon_{2k+1}})+\tilde{b}(G_{2k}e^{-\epsilon_{2k+1}}-\epsilon_{2k+1} [\tilde{a},\tilde{b}]F_{2k}-\epsilon_{2k+1} [\tilde{b},\tilde{b}]G_{2k})\\
&=\tilde{a}(F_{2k}e^{-\epsilon_{2k+1}})+\tilde{b}(G_{2k}(e^{-\epsilon_{2k+1}}-\epsilon_{2k+1} [\tilde{b},\tilde{b}])-\epsilon_{2k+1} [\tilde{a},\tilde{b}]F_{2k})\\
&=\tilde{a}(F_{2k}e^{-\epsilon_{2k+1}})+\tilde{b}(G_{2k}(e^{-\epsilon_{2k+1}}-\epsilon_{2k+1} \begin{pmatrix}
-z & 0\\
 0& z
 \end{pmatrix})-\epsilon_{2k+1} [\tilde{a},\tilde{b}]F_{2k})\\
&=\tilde{a}(F_{2k}e^{-\epsilon_{2k+1}})+\tilde{b}(G_{2k}e^{\epsilon_{2k+1}}-\epsilon_{2k+1} [\tilde{a},\tilde{b}]F_{2k})\\
\end{split}
\end{equation*}
which implies that (iii) holds for  $n=k$.
\end{proof}

\begin{lemma}\label{F-G}
%there exist $f_{2n}(M,\tilde{\lambda})$ and   $g_{2n}(M,\tilde{\lambda})$ such that 
 $$F_{2n}=\begin{pmatrix}
f_{2n}(M,\tilde{\lambda})  & 0\\
  0& f_{2n}(M^{-1},\tilde{\lambda}) 
 \end{pmatrix}, \quad -[\tilde{a},\tilde{b}]^{-1}G_{2n}=\begin{pmatrix}
g_{2n}(M,\tilde{\lambda})  & 0\\
  0& g_{2n}(M^{-1},\tilde{\lambda}) 
 \end{pmatrix}$$
where
$$f_{2n}(M,\tilde{\lambda})= \displaystyle\sum_{k=0}^{n} \tilde{c}_{2k}^{2n} \tilde{\lambda}^k,\quad  g_{2n}(M,\tilde{\lambda})= \displaystyle\sum_{k=0}^{n-1} \tilde{c}_{2k+1}^{2n} \tilde{\lambda}^k.$$
\end{lemma}
\begin{proof}
By (i) and (ii) of Proposition \ref{kmot-rep} and Lemma \ref{a-a-b},
\begin{equation*}
\begin{split}
F_{2}&=F_{1}e^{\epsilon_{2}}-\epsilon_{2}[\tilde{b}, \tilde{a}]G_{1}
=e^{-\epsilon_1}e^{\epsilon_{2}}+\epsilon_1\epsilon_{2}[\tilde{b}, \tilde{a}][\tilde{a}, \tilde{b}]=e^{-\epsilon_1+\epsilon_{2}}+\epsilon_1\epsilon_{2}\tilde{\lambda} I\\
& =\begin{pmatrix}
M^{-\epsilon_1+\epsilon_{2}}+\epsilon_1\epsilon_{2}\tilde{\lambda}  & 0\\
  0& M^{\epsilon_1-\epsilon_{2}}+\epsilon_1\epsilon_{2}\tilde{\lambda}
 \end{pmatrix}\\
&=\begin{pmatrix}
\tilde{c}_{0}^{2}(M,\epsilon)+\tilde{c}_{2}^{2}(M,\epsilon) \tilde{\lambda}  & 0\\
  0& \tilde{c}_{0}^{2}(M^{-1},\epsilon)+\tilde{c}_{2}^{2}(M^{-1},\epsilon) \tilde{\lambda} 
 \end{pmatrix}\\
&=\begin{pmatrix}
f_{2}(M,\tilde{\lambda})  & 0\\
  0& f_{2}(M^{-1},\tilde{\lambda}) 
 \end{pmatrix}\\
G_{2}&=G_{1}e^{-\epsilon_{2}}=-\epsilon_1[\tilde{a}, \tilde{b}]e^{-\epsilon_{2}}=-[\tilde{a}, \tilde{b}]\begin{pmatrix}
\tilde{c}_{1}^{2} (M,\epsilon)\tilde{\lambda}  & 0\\
  0& \tilde{c}_{1}^{2} (M^{-1},\epsilon) \tilde{\lambda}
 \end{pmatrix}\\
  &=-[\tilde{a}, \tilde{b}]\begin{pmatrix}
g_{2}(M,\tilde{\lambda})  & 0\\
  0& g_{2}(M^{-1},\tilde{\lambda}) 
 \end{pmatrix},\\
\end{split}
\end{equation*}
which implies that the statement holds for $n=1$. Assume that the statement holds for all $n\leq m$ to proceed by induction on $n$. Then 
\begin{equation*}
\begin{split}
F_{2m+2}&=F_{2m+1}e^{\epsilon_{2m+2}}-\epsilon_{2m+2}[\tilde{b}, \tilde{a}]G_{2m+1} \\
&=F_{2m}e^{-\epsilon_{2m+1}}e^{\epsilon_{2m+2}}-\epsilon_{2m+2}[\tilde{b}, \tilde{a}]
(G_{2m}e^{\epsilon_{2m+1}}-\epsilon_{2m+1}[\tilde{a}, \tilde{b}]F_{2m}) \\
&=F_{2m}(e^{-\epsilon_{2m+1}+\epsilon_{2m+2}}+\epsilon_{2m+1}\epsilon_{2m+2}\tilde{\lambda} I)-\epsilon_{2m+2}G_{2m}e^{\epsilon_{2m+1}}[\tilde{b}, \tilde{a}]\\
&=\begin{pmatrix}
f_{2m+2}(M,\tilde{\lambda})  & 0\\
  0& f_{2m+2}(M^{-1},\tilde{\lambda}) 
 \end{pmatrix}\\
\end{split}
\end{equation*}
for some $f_{2m+2}(M,\tilde{\lambda}) $ and  thus \begin{equation*}
\begin{split}
f_{2m+2}&=
\displaystyle\sum_{k=0}^{m} \tilde{c}_{2k}^{2m}M^{-\epsilon_{2m+1}+\epsilon_{2m+2}}\tilde{\lambda}^k+\epsilon_{2m+1}\epsilon_{2m+2}\sum_{k=0}^{m} \tilde{c}_{2k}^{2m}\tilde{\lambda}^{k+1}+\epsilon_{2m+2}M^{\epsilon_{2m+1}}\sum_{k=0}^{m-1} \tilde{c}_{2k+1}^{2m}\tilde{\lambda}^{k+1} \\
&=
\tilde{c}_{0}^{2m+2}+\displaystyle\sum_{k=1}^{m}(\tilde{c}_{2k}^{2m}M^{-\epsilon_{2m+1}+\epsilon_{2m+2}}+\epsilon_{2m+1}\epsilon_{2m+2}\tilde{c}_{2k-2}^{2m}+\epsilon_{2m+2}M^{\epsilon_{2m+1}}\tilde{c}_{2k-1}^{2m})\tilde{\lambda}^{k}\\
&=
\tilde{c}_{0}^{2m+2}+\displaystyle\sum_{k=1}^{m}\tilde{c}_{2k}^{2m+2}\tilde{\lambda}^{k} \\
&=\displaystyle\sum_{k=0}^{m}\tilde{c}_{2k}^{2m+2}\tilde{\lambda}^{k} \\
\end{split}
\end{equation*}
which implies  that $F_{2m+2}$ also satisfies the formula. Note that the equality
$$\tilde{c}_{2k}^{2m}M^{-\epsilon_{2m+1}+\epsilon_{2m+2}}+\epsilon_{2m+1}\epsilon_{2m+2}\tilde{c}_{2k-2}^{2n}+\epsilon_{2m+2}M^{\epsilon_{2m+1}}\tilde{c}_{2k-1}^{2m}=\tilde{c}_{2k}^{2m+2}$$
is (v) of Lemma \ref{lem-2}.

The formula for $G_{2n}$  is also proved similarly using Lemma \ref{lem-2}, (vi) as follows.  Since 
\begin{equation*}
\begin{split}
-[\tilde{a}, \tilde{b}]^{-1}&G_{2m+2}=
-[\tilde{a}, \tilde{b}]^{-1}G_{2m+1}e^{-\epsilon_{2m+2}}\\
&=-[\tilde{a}, \tilde{b}]^{-1}(G_{2m}e^{\epsilon_{2m+1}}-\epsilon_{2m+1}[\tilde{a}, \tilde{b}]F_{2m})e^{-\epsilon_{2m+2}}\\
&=-[\tilde{a}, \tilde{b}]G_{2m}e^{\epsilon_{2m+1}-\epsilon_{2m+2}}+\epsilon_{2m+1}F_{2m}e^{-\epsilon_{2m+2}}\\
&=\begin{pmatrix}
g_{2m}(M,\tilde{\lambda})  & 0\\
  0& g_{2m}(M^{-1},\tilde{\lambda}) 
 \end{pmatrix}
e^{\epsilon_{2m+1}-\epsilon_{2m+2}}
+\epsilon_{2m+1}
\begin{pmatrix}
f_{2m}(M,\tilde{\lambda})  & 0\\
  0& f_{2m}(M^{-1},\tilde{\lambda}) 
 \end{pmatrix}e^{-\epsilon_{2m+2}}\\
\end{split}
\end{equation*}
we have 
\begin{equation*}
\begin{split}
g_{2m+2}&=g_{2m}M^{\epsilon_{2m+1}-\epsilon_{2m+2}}+\epsilon_{2m+1}f_{2m}M^{-\epsilon_{2m+2}}\\
&=\displaystyle\sum_{k=0}^{m-1}\tilde{c}_{2k+1}^{2m}M^{\epsilon_{2m+1}-\epsilon_{2m+2}}\tilde{\lambda}^{k}+\epsilon_{2m+1}\displaystyle\sum_{k=0}^{m}\tilde{c}_{2k}^{2m}M^{-\epsilon_{2m+2}}\tilde{\lambda}^{k}\\
&=\displaystyle\sum_{k=0}^{m-1}(\tilde{c}_{2k+1}^{2m}M^{\epsilon_{2m+1}-\epsilon_{2m+2}}+\epsilon_{2m+1}\tilde{c}_{2k}^{2m}M^{-\epsilon_{2m+2}})\tilde{\lambda}^{k}+\epsilon_{2m+1}\tilde{c}_{2m}^{2m}M^{-\epsilon_{2m+2}}\tilde{\lambda}^{m}
\\
&=\displaystyle\sum_{k=0}^{m-1}(\tilde{c}_{2k+1}^{2m}M^{\epsilon_{2m+1}-\epsilon_{2m+2}}+\epsilon_{2m+1}\tilde{c}_{2k}^{2m}M^{-\epsilon_{2m+2}})\tilde{\lambda}^{k}+\tilde{c}_{2m}^{2m+2}\tilde{\lambda}^{m}\\
&=\displaystyle\sum_{k=0}^{m-1}\tilde{c}_{2k+1}^{2m+2}\tilde{\lambda}^{k}+\tilde{c}_{2m}^{2m+2}\tilde{\lambda}^{m}\\
&=\displaystyle\sum_{k=0}^{m}\tilde{c}_{2k+1}^{2m+2}\tilde{\lambda}^{k},
\end{split}
\end{equation*}
which implies  that $G_{2m+2}$ also satisfies the formula. 
\end{proof}
\begin{remark}
It is quite remarkable to notice that the formulas for $f_{\alpha-1}$ and $g_{\alpha-1}$ in Lemma \ref{F-G} are obtained from the formulas $W_{11}$ and $W_{12}$ in Proposition \ref{W-eq}  replacing $\lambda$ and $\widehat{i_1 \cdots i_{k}}$ by $\tilde{\lambda}$ and $\widehat{i_1 \cdots i_{k}}^{alt}$ respectively.
That is,
\begin{equation*}
\begin{split}
W_{11}&=\sum_{k=0}^{\frac{\alpha-1}{2}} c_{2k}^{\alpha-1} \lambda^k\\
&= M^{\sum_i \epsilon_i} +\lambda\sum_{i_1:\text {odd}}^{\wedge} \epsilon_{i_1} \epsilon_{i_2}M^{\widehat{i_1 i_2}}+ \cdots+\lambda^k\sum_{i_1:\text {odd}}^{\wedge} \epsilon_{i_1} \cdots\epsilon_{i_{2k}}M^{\widehat{i_1 \cdots i_{2k}}}+ \cdots +\lambda^{\frac{\alpha-1}{2}} \epsilon_{1} \epsilon_{2}\epsilon_{3}\cdots\epsilon_{\alpha-1} \\
W_{12}&=\sum_{k=0}^{\frac{\alpha-3}{2}} c_{2k+1}^{\alpha-1} \lambda^k=\sum_{k=0}^{\frac{\alpha-3}{2}} d_{2k+1}^{\alpha-1} \lambda^k\\
&= \sum_{i_1:\text {odd}}^{\wedge} \epsilon_{i_1} M^{\widehat{i_1}}+ \lambda\sum_{i_1:\text {odd}}^{\wedge} \epsilon_{i_1} \epsilon_{i_2}\epsilon_{i_3}M^{\widehat{i_1 i_2 i_3 }}+ \cdots+\lambda^{\frac{\alpha-3}{2}}\epsilon_{1} \epsilon_{2}\epsilon_{3}\cdots\epsilon_{\alpha-2} M^{-\epsilon_{\alpha-1}}    \\
\end{split}
\end{equation*}
and 
\begin{equation*}
\begin{split}
f_{\alpha-1}&(M,\tilde{\lambda})=\displaystyle\sum_{k=0}^{\frac{\alpha-1}{2}} \tilde{c}_{2k}^{\alpha-1} \tilde{\lambda}^k\\
&=1 +\tilde{\lambda} \sum_{i_1:\text {odd}}^{\wedge} \epsilon_{i_1} \epsilon_{i_2}M^{\widehat{i_1 i_2}^{alt}}+ \cdots+\tilde{\lambda} ^k\sum_{i_1:\text {odd}}^{\wedge} \epsilon_{i_1} \cdots\epsilon_{i_{2k}}M^{\widehat{i_1 \cdots i_{2k}}^{alt}}+ \cdots +\tilde{\lambda} ^{\frac{\alpha-1}{2}} \epsilon_{1} \epsilon_{2}\epsilon_{3}\cdots\epsilon_{\alpha-1} \\
g_{\alpha-1}&(M,\tilde{\lambda})=\displaystyle\sum_{k=0}^{\frac{\alpha-3}{2}} \tilde{c}_{2k+1}^{\alpha-1}\tilde{\lambda}^k=\displaystyle\sum_{k=0}^{\frac{\alpha-3}{2}} \tilde{d}_{2k+1}^{\alpha-1}\tilde{\lambda}^k\\
&=\sum_{i_1:\text {odd}}^{\wedge} \epsilon_{i_1} M^{\widehat{i_1}^{alt}}+ \tilde{\lambda}\sum_{i_1:\text {odd}}^{\wedge} \epsilon_{i_1} \epsilon_{i_2}\epsilon_{i_3}M^{\widehat{i_1 i_2 i_3 }^{alt}}+ \cdots+\tilde{\lambda}^{\frac{\alpha-3}{2}}\epsilon_{1} \epsilon_{2}\epsilon_{3}\cdots\epsilon_{\alpha-2} M^{-\epsilon_{\alpha-1}}.
\end{split}
\end{equation*}
Note also that each root of $f_{\alpha-1}(1,\tilde{\lambda})$ corresponds to a parabolic representation, and 
both $f_{\alpha-1}$ and $\sqrt{-\lambda}g_{\alpha-1}$ are $\epsilon$-Chebyshev polynomials when $M=1$ by Theorem 5.1 of  \cite{Jo-Kim2}.
\end{remark}
%From Lemma \ref{F-G} and Corollary \ref{F-f}.
\begin{corollary}\label{F-f}
\begin{enumerate}
\item [\rm (i)]
$f_{\alpha-1}(M^{-1},\tilde{\lambda})=f_{\alpha-1}(M,\tilde{\lambda})$
\item [\rm (ii)]  $f_{\alpha-1}(M,\tilde{\lambda})\in \mathbb Z[z^2,\tilde{\lambda}]= \mathbb Z[\lambda, \tilde{\lambda}] $ 
\item [\rm (iii)]
 $F_{\alpha-1}=\begin{pmatrix}
f_{\alpha-1}(M,\tilde{\lambda})  & 0\\
  0& f_{\alpha-1}(M,\tilde{\lambda}) 
 \end{pmatrix}$
\item [\rm (iv)] $G_{\alpha-1}=-\begin{pmatrix}
g_{\alpha-1}(M,\tilde{\lambda})  & 0\\
  0& g_{\alpha-1}(M^{-1},\tilde{\lambda}) 
 \end{pmatrix}[\tilde{a}, \tilde{b}]$
\end{enumerate}
where
\begin{equation*}
f_{\alpha-1}(M,\tilde{\lambda})=\displaystyle\sum_{k=0}^{\frac{\alpha-1}{2}} \tilde{c}_{2k}^{\alpha-1} \tilde{\lambda}^k,\quad 
g_{\alpha-1}(M,\tilde{\lambda})=\displaystyle\sum_{k=0}^{\frac{\alpha-3}{2}} \tilde{c}_{2k+1}^{\alpha-1}\tilde{\lambda}^k=\displaystyle\sum_{k=0}^{\frac{\alpha-3}{2}} \tilde{d}_{2k+1}^{\alpha-1}\tilde{\lambda}^k.
\end{equation*}
\end{corollary}
\begin{proof}
The statement (iv) is obvious from Lemma \ref{F-G}.
By Lemma \ref{symmetric}, (iii),
$$f_{\alpha-1}(M^{-1},\tilde{\lambda})=
\displaystyle\sum_{k=0}^{m}\tilde{c}_{2k}^{2m}(M^{-1},\tilde{\lambda})\tilde{\lambda}^{k}
=\displaystyle\sum_{k=0}^{m}\tilde{c}_{2k}^{2m}(M,\tilde{\lambda})\tilde{\lambda}^{k}
=f_{\alpha-1}(M,\tilde{\lambda}),$$
which proves (i)  and (iii).
(ii) $f_{\alpha-1}(M,\tilde{\lambda}) \in \mathbb Z[z^2,\tilde{\lambda}]=\mathbb Z[\lambda,\tilde{\lambda}]$ follows from  $z^2=\tilde{\lambda}-\lambda$.
\end{proof}
\begin{corollary}\label{G-det} Let $\bar{g}_{\alpha-1}(M,\tilde{\lambda}):=g_{\alpha-1}(M^{-1},\tilde{\lambda})$. Then 
\begin{equation*}
z(1-f_{\alpha-1}^2+\tilde{\lambda}g_{\alpha-1}\bar{g}_{\alpha-1}) =\tilde{\lambda}f_{\alpha-1}( g_{\alpha-1}-\bar{g}_{\alpha-1}).
\end{equation*}
Especially, 
if $f_{\alpha-1}(M,\tilde{\lambda})=0$, then 
\begin{equation}\label{g-alpha}
\det G_{\alpha-1}=-\tilde{\lambda}g_{\alpha-1}(M,\tilde{\lambda})g_{\alpha-1}(M^{-1},\tilde{\lambda})=1.
\end{equation}
\end{corollary}
\begin{proof}
From Proposition \ref{kmot-rep}, we have $a, b \in SL(2,\mathbb C)$ such that 
\begin{equation*}
\begin{split}
\tilde{c}_{\alpha-1}\sim \tilde{a}\begin{pmatrix}f_{\alpha-1}(M,\tilde{\lambda})  & 0\\
  0& f_{\alpha-1}(M,\tilde{\lambda}) 
 \end{pmatrix}-\tilde{b}
\begin{pmatrix}g_{\alpha-1}(M,\tilde{\lambda})  & 0\\
  0& g_{\alpha-1}(M^{-1},\tilde{\lambda}) 
 \end{pmatrix}
[\tilde{a},\tilde{b}].
\end{split}
\end{equation*}
Thus 
\begin{equation*}
\begin{split}
&\tilde{c}_{\alpha-1}\sim (c_1,c_3)\\
&c_1= f_{\alpha-1}a_1-\langle a_1,b_3\rangle g_{\alpha-1}b_1,\,\, 
c_3= f_{\alpha-1}a_3-\langle a_3,b_1\rangle \bar{g}_{\alpha-1}b_3.
\end{split}
\end{equation*}
Hence we have
\begin{equation*}
\begin{split}
-z&=\langle c_1,c_3\rangle=\langle  f_{\alpha-1}a_1-\langle a_1,b_3\rangle g_{\alpha-1}b_1, f_{\alpha-1}a_3-\langle a_3,b_1\rangle \bar{g}_{\alpha-1}b_3\rangle\\
&=f_{\alpha-1}^2\langle a_1,a_3\rangle+\langle a_1,b_3\rangle \langle a_3,b_1\rangle (g_{\alpha-1}\bar{g}_{\alpha-1}\langle b_1,b_3\rangle -f_{\alpha-1}\bar{g}_{\alpha-1} +f_{\alpha-1} g_{\alpha-1})\\
&=-zf_{\alpha-1}^2+\tilde{\lambda}(zg_{\alpha-1}\bar{g}_{\alpha-1} +f_{\alpha-1}\bar{g}_{\alpha-1} -f_{\alpha-1} g_{\alpha-1}),
\end{split}
\end{equation*}
which implies that
\begin{equation*}
z(1-f_{\alpha-1}^2+\tilde{\lambda}g_{\alpha-1}\bar{g}_{\alpha-1}) =\tilde{\lambda}f_{\alpha-1}( g_{\alpha-1}-\bar{g}_{\alpha-1}).
\end{equation*}
This completes the proof.
\end{proof}
Since we can normalize any two matrices $A, B$ in 
$SL(2,\mathbb C)$, which have the same trace and do not commute, into a particular pair $$ m=\begin{pmatrix}
M & 1\\
  0& \frac{1}{M} 
 \end{pmatrix},\quad \begin{pmatrix}
M & 0\\
  \lambda& \frac{1}{M} 
 \end{pmatrix},\quad \lambda=2-tr(AB^{-1}),$$
by Lemma \ref{Riley-lemma-3}, we will assume from now on  \begin{equation*}
\rho(x)=\begin{pmatrix}
M & 1\\
  0 & \frac{1}{M} 
 \end{pmatrix} \quad \text{and}\quad
\rho(y)=\begin{pmatrix}
M & 0\\
  \lambda& \frac{1}{M} 
 \end{pmatrix}
\end{equation*}
and 
$a= \begin{pmatrix}
1 & 0\\
 0& 1
 \end{pmatrix}$ and 
$b= \begin{pmatrix}
z & 1\\
\lambda& \frac{\lambda+1}{z}
 \end{pmatrix}$ when $\rho : F_{x,y}  \rightarrow SL(2,\bc)$ is non-abelian.
So we have 
\begin{equation*}
\begin{split}
\tilde{c}_{\alpha-1}\sim \tilde{a}\begin{pmatrix}f_{\alpha-1}(M,\tilde{\lambda})  & 0\\
  0& f_{\alpha-1}(M,\tilde{\lambda}) 
 \end{pmatrix}-\tilde{b}
\begin{pmatrix}g_{\alpha-1}(M,\tilde{\lambda})  & 0\\
  0& g_{\alpha-1}(M^{-1},\tilde{\lambda}) 
 \end{pmatrix}
\begin{pmatrix}
-1  & 0\\
  0& \tilde{\lambda}
\end{pmatrix}
\end{split}
\end{equation*}
by Proposition \ref{kmot-rep} and Lemma \ref{F-G}.
\begin{theorem}\label{kmot-rep2}
Let $\rho : F_{x,y}  \rightarrow SL(2,\bc)$ be a group homomorphism satisfying 
\begin{equation*}
\rho(x)=\begin{pmatrix}
M & 1\\
  0 & \frac{1}{M} 
 \end{pmatrix} \quad \text{and}\quad
\rho(y)=\begin{pmatrix}
M & 0\\
  \lambda& \frac{1}{M} 
 \end{pmatrix}.
\end{equation*} Then  
$\rho$ defines a non-abelian representation of  a kmot group $G(\epsilon)$ if and only if  $f_{\alpha-1}(M,\tilde{\lambda})=0$. Furthermore, $f_{\alpha-1}(M,\tilde{\lambda})$ is the Riley polynomial of $G(\epsilon)$.
\end{theorem}
\begin{proof}
$\rho$ defines a non-abelian representation of  a kmot group $G$ if and only if $\tilde{c}_{\alpha-1}\sim  \tilde{b}$, since 
$c_{\alpha-1}mc_{\alpha-1}^{-1}=\rho(x_{\alpha-1})=\rho(wxw^{-1})=\rho(y)=bmb^{-1}.$
Since
 $$\tilde{c}_{\alpha-1}\sim \tilde{a}F_{\alpha-1}+\tilde{b}G_{\alpha-1}=f_{\alpha-1}(M,\tilde{\lambda})\begin{pmatrix}
1 & 1\\
 0& -z
 \end{pmatrix}+\begin{pmatrix}
z & 0\\
\lambda& -1
 \end{pmatrix}\begin{pmatrix}
g_{\alpha-1}(M,\tilde{\lambda})  & 0\\
  0& -\tilde{\lambda}g_{\alpha-1}(M^{-1},\tilde{\lambda}) 
 \end{pmatrix}$$
by Proposition \ref{kmot-rep}, 
$\tilde{c}_{\alpha-1}\sim  \tilde{b}$ is equivalent to  
\begin{equation}\label{c-b}
f_{\alpha-1}(M,\tilde{\lambda})\begin{pmatrix}
1 & 1\\
 0& -z
 \end{pmatrix}+\begin{pmatrix}
z & 0\\
\lambda& -1
 \end{pmatrix}\begin{pmatrix}
g_{\alpha-1}(M,\tilde{\lambda})  & 0\\
  0& -\tilde{\lambda}g_{\alpha-1}(M^{-1},\tilde{\lambda}) 
 \end{pmatrix}\sim \begin{pmatrix}
z & 0\\
\lambda& -1
 \end{pmatrix},
\end{equation}
which holds when $f_{\alpha-1}(M,\tilde{\lambda})=0$ and $-\tilde{\lambda}g_{\alpha-1}(M,\tilde{\lambda})g_{\alpha-1}(M^{-1},\tilde{\lambda})=1$.
But  if $f_{\alpha-1}(M,\tilde{\lambda})=0$ then $-\tilde{\lambda} g_{\alpha-1}(M,\tilde{\lambda})g_{\alpha-1}(M^{-1},\tilde{\lambda})=1$ is automatically satisfied by Corollary \ref{G-det}, which implies that 
$$f_{\alpha-1}(M,\tilde{\lambda})\,|\, 1+\tilde{\lambda}g_{\alpha-1}(M,\tilde{\lambda})g_{\alpha-1}(M^{-1},\tilde{\lambda})$$
and thus $\rho$ defines a non-abelian representation of  a kmot group $G$ with (\ref{kmot-group}) and (\ref{presentation1}) if and only if $f_{\alpha-1}(M,\tilde{\lambda})=0$. 
Since the $\lambda$-degrees of $W_{11}-zW_{12}$ and $f_{\alpha-1}(M,\tilde{\lambda})$ are both $\frac{\alpha-1}{2}$ and $\tilde{c}_{\alpha-1}^{\alpha-1}=c_{\alpha-1}^{\alpha-1}=\epsilon_{1} \epsilon_{2}\epsilon_{3}\cdots\epsilon_{\alpha-1} $, they must be the same by  Proposition \ref{Riley poly}.  Therefore $f_{\alpha-1}(M,\tilde{\lambda})$ is the Riely polynomial of $G(\epsilon)$.
%But we can check that $W_{11}-zW_{12}=f_{\alpha-1}(M,\tilde{\lambda})$ by comparing their degrees, so $f_{\alpha-1}(M,\tilde{\lambda})$ divides $1+\tilde{\lambda}g_{\alpha-1}(M,\tilde{\lambda})g_{\alpha-1}(M^{-1},\tilde{\lambda})$, which completes the proof.
\end{proof}
\begin{remark}
Since $f_{\alpha-1}(M^{-1},\tilde{\lambda})=f_{\alpha-1}(M,\tilde{\lambda})$ by Corollary \ref{F-f}, if $(M,\tilde{\lambda})$ is a root of $f_{\alpha-1}$ then so is $(M^{-1},\tilde{\lambda})$ and they give the same  non-abelian representation of  a kmot group $G(\epsilon)$ up to conjugation.
But  $f_{\alpha-1}\in \mathbb Z[s,\tilde{\lambda}], s:=z^2$ and  thus there is a 1-1 correspondence between the set of roots $\{(s,\tilde{\lambda})\}$ of $f_{\alpha-1}$ and the set of non-abelian representation of  a kmot group $G(\epsilon)$ up to conjugation.  
\end{remark}

\begin{proposition}\label{w-tilde}
$   W_{11}=f_{\alpha-1}(M,\tilde{\lambda})+zg_{\alpha-1}(M,\tilde{\lambda}),\quad  W_{12}=g_{\alpha-1}(M,\tilde{\lambda})$
\end{proposition} 
\begin{proof}
Since $c_{\alpha-1}mc_{\alpha-1}^{-1}=\rho(x_{\alpha-1})=\rho(wxw^{-1})=W\rho(x)W^{-1}$, we have 
$\tilde{W}\sim \tilde{c}_{\alpha-1},$ that is 
\begin{equation*}
\begin{split}
\begin{pmatrix}
W_{11} & W_{11}-zW_{12}\\
\lambda W_{12}& \lambda W_{12}-zW_{22}
 \end{pmatrix}&\sim 
f_{\alpha-1}(M,\tilde{\lambda})\begin{pmatrix}
1 & 1\\
 0& -z
 \end{pmatrix}+\begin{pmatrix}
z & 0\\
\lambda& -1
 \end{pmatrix}\begin{pmatrix}
g_{\alpha-1}(M,\tilde{\lambda})  & 0\\
  0& -\tilde{\lambda}g_{\alpha-1}(M^{-1},\tilde{\lambda}) 
 \end{pmatrix}\\
&=\begin{pmatrix}
f_{\alpha-1}(M,\tilde{\lambda})+zg_{\alpha-1}(M,\tilde{\lambda})  & f_{\alpha-1}(M,\tilde{\lambda})\\
  \lambda g_{\alpha-1}(M,\tilde{\lambda})& -zf_{\alpha-1}(M,\tilde{\lambda})+\tilde{\lambda}g_{\alpha-1}(M^{-1},\tilde{\lambda}) 
 \end{pmatrix}
\end{split}
\end{equation*}
Hence $W_{12}=g_{\alpha-1}(M,\tilde{\lambda})$ if and only if  $W_{11}-zW_{12}=f_{\alpha-1}(M,\tilde{\lambda})$. But since
$W_{11}-zW_{12}=f_{\alpha-1}(M,\tilde{\lambda})$ has been already proved in Theorem \ref{kmot-rep2},  $W_{12}=g_{\alpha-1}(M,\tilde{\lambda})$ and thus $   W_{11}=f_{\alpha-1}(M,\tilde{\lambda})+zg_{\alpha-1}(M,\tilde{\lambda})$ are also proved.
\end{proof}

\begin{remark}
Note that Theorem \ref{kmot-rep2} and 
Proposition \ref{w-tilde}
follow immediately from Proposition \ref{prop-1}, Proposition \ref{W-eq}, and Proposition \ref{Riley poly}.
\end{remark}

\section{Recursive formulas for Riley polynomials and Alexander polynomials}%=========================
%\textcolor{blue}{From Proposition \ref{kmot-rep} and Theorem \ref{kmot-rep2}, 
In the previous section, 
we have seen that there are two polynomials  $f_{\alpha-1}\in \mathbb Z[\lambda, \tilde{\lambda}] $ and $g_{\alpha-1}\in \mathbb Z[M, M^{-1}, \tilde{\lambda}] $  for a given symmetric $\epsilon_i$-sequence of length $\alpha-1$. Since these polynomials are different in general for distinct  $\epsilon_i$-sequences, let us use the notation
 $f_{\epsilon}$ and $g_{\epsilon}$ instead. Using this notation we have 
\begin{equation}\label{f-g}
f_{-\epsilon}(M,\tilde{\lambda})=f_{\epsilon}(M,\tilde{\lambda})=f_{\epsilon}(M^{-1},\tilde{\lambda}), \quad g_{-\epsilon}(M,\tilde{\lambda})= -g_{\epsilon}(M^{-1},\tilde{\lambda}). 
\end{equation}
for any symmetric $\epsilon_i$-sequence ${\bf \epsilon}=(\epsilon_1,\epsilon_2,\cdots,\epsilon_{\alpha-1})$.
\begin{proposition}\label{w-quandle}
Let ${\bf \epsilon}, {\bf \epsilon'}, {\bf \epsilon''}$ be symmetric $\epsilon_i$-sequences such that
$${\bf \epsilon}=(\epsilon_1,\epsilon_2,\cdots,\epsilon_{\alpha-1}), \quad  {\bf \epsilon'}=(\epsilon_2,\epsilon_3,\cdots,\epsilon_{\alpha-2}), \quad {\bf \epsilon''}=(\epsilon_3,\epsilon_4,\cdots,\epsilon_{\alpha-3}).$$ Then 
\begin{enumerate}
\item [\rm (i)]
$g_{\epsilon}=  g_{\epsilon'}(M^{-1},\tilde{\lambda})M^{-2\epsilon_1}+\epsilon_1f_{\epsilon'}M^{-\epsilon_1}$
\item [\rm (ii)]
$W(\epsilon)=\begin{pmatrix}
f_{\epsilon}+zg_{\epsilon} & g_{\epsilon}\\
  \lambda g_{\epsilon} & f_{\epsilon'}-zg_{\epsilon}
 \end{pmatrix}$
\item [\rm (iii)]
 $\tilde{\lambda}(g_{\epsilon}(M^{-1},\tilde{\lambda})-g_{\epsilon}(M,\tilde{\lambda}))=z(f_{\epsilon}(\lambda, \tilde{\lambda})-f_{\epsilon'}(\lambda, \tilde{\lambda}))$
\item [\rm (iv)] $1-\det G_{\epsilon}=1+\tilde{\lambda}g_{\epsilon}(M,\tilde{\lambda})g_{\epsilon}(M^{-1},\tilde{\lambda})=f_{\epsilon}(\lambda, \tilde{\lambda})f_{\epsilon'}(\lambda, \tilde{\lambda})$
\item [\rm (v)]
$f_{\epsilon}= (f_{\epsilon''}M^{2\epsilon_1}+\lambda f_{\epsilon'})+
\bar{g}_{\epsilon'}(2\epsilon_1\lambda M^{\epsilon_1}+z(-\lambda+M^{2\epsilon_1}) )-zg_{\epsilon}$  for $\alpha>3$.
\end{enumerate}
Here $\bar{g}_{\epsilon}(M,\tilde{\lambda}):=g_{\epsilon}(M^{-1},\tilde{\lambda})$. We let $f_{\epsilon'}(\lambda, \tilde{\lambda})=1$ and 
$g_{\epsilon'}(M, \tilde{\lambda})=0$
 when ${\bf \epsilon}=(\epsilon_1,\epsilon_1)$.
\end{proposition}
\begin{proof}
We use induction on $\alpha$.
When $\alpha=3$, that is, ${\bf \epsilon}=(\epsilon_1,\epsilon_1)$,  we have 
$$W(\epsilon)=\begin{pmatrix}
M^{2\epsilon_1}+\lambda & \epsilon_1M^{-\epsilon_1}\\
  \epsilon_1\lambda M^{-\epsilon_1}& M^{-2\epsilon_1} 
 \end{pmatrix}.$$
Hence we can check that  statements are true for $\alpha=3$ as follows.
\begin{equation*}
\begin{split}
W_{12}(\epsilon)&=\epsilon_1M^{-\epsilon_1}=g_{\epsilon}(M,\tilde{\lambda})\\ 
W_{11}(\epsilon)-zW_{12}(\epsilon)&=M^{2\epsilon_1}+\lambda -z\epsilon_1M^{-\epsilon_1}=
M^{2\epsilon_1}+\lambda  -(M^{\epsilon_1}-M^{-\epsilon_1})M^{-\epsilon_1}\\
&=M^{2\epsilon_1}+M^{-2\epsilon_1}-1+\lambda =M^2+M^{-2}-1+\lambda\\
&=(z^2+2)-1+\lambda=1+\tilde{\lambda}=f_{\epsilon}(M,\tilde{\lambda})\\
W_{22}(\epsilon)&=M^{-2\epsilon_1}=f_{\epsilon'}(\lambda, \tilde{\lambda}) \\
g_{\epsilon}(M,\tilde{\lambda})&=\epsilon_1M^{-\epsilon_1}=g_{\epsilon'}(M^{-1},\tilde{\lambda})M^{-2\epsilon_1}+\epsilon_1f_{\epsilon'}M^{-\epsilon_1}\\
1+\tilde{\lambda}g_{\epsilon}(M,\tilde{\lambda})g_{\epsilon}(M^{-1},\tilde{\lambda})&=1+\tilde{\lambda}\epsilon_1M^{-\epsilon_1}\epsilon_1M^{\epsilon_1}=1+\tilde{\lambda}=f_{\epsilon}(\lambda, \tilde{\lambda})f_{\epsilon'}(\lambda, \tilde{\lambda})
\end{split}
\end{equation*}
Now it suffices to show that the statements are true when they hold for $\epsilon'$.

Since
\begin{equation*}
\begin{split}
W(\epsilon)&= \begin{pmatrix}
M^{\epsilon_1} & \epsilon_1\\
  0 & M^{-\epsilon_1}
 \end{pmatrix}W^*(\epsilon') 
\begin{pmatrix}
M^{\epsilon_1} &0\\
  \epsilon_1\lambda &M^{-\epsilon_1} 
 \end{pmatrix}\\
&= \begin{pmatrix}
M^{\epsilon_1} & \epsilon_1\\
 0 &M^{-\epsilon_1} 
 \end{pmatrix}
 \begin{pmatrix}
 f_{\epsilon''}+z\bar{g}_{\epsilon'}& \bar{g}_{\epsilon'}\\
  \lambda\bar{g}_{\epsilon'}& f_{\epsilon'}-z\bar{g}_{\epsilon'}
 \end{pmatrix}
\begin{pmatrix}
M^{\epsilon_1} &0\\
  \epsilon_1\lambda &M^{-\epsilon_1} 
 \end{pmatrix}\\
&=\begin{pmatrix}
*& \bar{g}_{\epsilon'}M^{-2\epsilon_1}+\epsilon_1f_{\epsilon'}M^{-\epsilon_1}\\
  \lambda(\bar{g}_{\epsilon'}+\epsilon_1(f_{\epsilon'}-z\bar{g}_{\epsilon'})M^{-\epsilon_1})&( f_{\epsilon'}-z\bar{g}_{\epsilon'})M^{-2\epsilon_1}\\
 \end{pmatrix}
 \end{split}
\end{equation*}
where
\begin{equation*}
*= (f_{\epsilon''}+z\bar{g}_{\epsilon'})M^{2\epsilon_1}+2\epsilon_1\lambda\bar{g}_{\epsilon'}M^{\epsilon_1} +\lambda(f_{\epsilon'}-z\bar{g}_{\epsilon'}),
\end{equation*}
  we have by Proposition \ref{w-tilde} that 
$$f_{\epsilon}+zg_{\epsilon}= (f_{\epsilon''}M^{2\epsilon_1}+\lambda f_{\epsilon'})+
\bar{g}_{\epsilon'}(2\epsilon_1\lambda M^{\epsilon_1}+z(-\lambda+M^{2\epsilon_1}) )$$
%\textcolor{blue}{
%\begin{equation*}
%\begin{split}
%f_{\epsilon}+zg_{\epsilon}&= (f_{\epsilon''}+z\bar{g}_{\epsilon'})M^{2\epsilon_1}+2\epsilon_1\lambda\bar{g}_{\epsilon'}M^{\epsilon_1} +\lambda(f_{\epsilon'}+z\bar{g}_{\epsilon'})\\
%&= (f_{\epsilon''}M^{2\epsilon_1}+\lambda f_{\epsilon'})+
%\bar{g}_{\epsilon'}(2\epsilon_1\lambda M^{\epsilon_1}+z(M^{2\epsilon_1}\lambda) )\\
% \end{split}
%\end{equation*}
%}
and
 $$g_{\epsilon}= \bar{g}_{\epsilon'}M^{-2\epsilon_1}+\epsilon_1f_{\epsilon'}M^{-\epsilon_1},$$  which prove (v) and (i).

To prove (ii), we just need to show  that $f_{\epsilon'}(\lambda,\tilde{\lambda})-zg_{\epsilon}(M,\tilde{\lambda})=( f_{\epsilon'}-z\bar{g}_{\epsilon'})M^{-2\epsilon_1}$.
\begin{equation*}
\begin{split}
f_{\epsilon'}(\lambda,\tilde{\lambda})-zg_{\epsilon}(M,\tilde{\lambda})&=f_{\epsilon'}-z(\bar{g}_{\epsilon'}M^{-2\epsilon_1}+\epsilon_1f_{\epsilon'}M^{-\epsilon_1})\\
%&=f_{\epsilon'}-\epsilon_1zM^{-\epsilon_1}f_{\epsilon'}-z\bar{g}_{\epsilon'}M^{-2\epsilon_1}\\
&=f_{\epsilon'}-(M^{\epsilon_1}-M^{-\epsilon_1})M^{-\epsilon_1}f_{\epsilon'}-z\bar{g}_{\epsilon'}M^{-2\epsilon_1}\\
&=M^{-2\epsilon_1}(f_{\epsilon'}-z\bar{g}_{\epsilon'})
 \end{split}
\end{equation*}

Since
\begin{equation}\label{det-w}
\begin{split}
1=\det W(\epsilon)&=
(f_{\epsilon}+zg_{\epsilon})( f_{\epsilon'}-zg_{\epsilon})-  \lambda g_{\epsilon}^2\\
&=f_{\epsilon} f_{\epsilon'}+zg_{\epsilon}(f_{\epsilon'}-f_{\epsilon})-z^2g_{\epsilon}^2-  \lambda g_{\epsilon}^2\\
&=f_{\epsilon} f_{\epsilon'}+zg_{\epsilon}(f_{\epsilon'}-f_{\epsilon})-\tilde{\lambda} g_{\epsilon}^2
 \end{split}
\end{equation}
we have
\begin{equation*}
\begin{split}
0&=\det W(\epsilon)-\det \tilde{W}(\epsilon)\\
&=(f_{\epsilon} f_{\epsilon'}+zg_{\epsilon}(f_{\epsilon'}-f_{\epsilon})-\tilde{\lambda} g_{\epsilon}^2)-(f_{\epsilon} f_{\epsilon'}-z\bar{g}_{\epsilon}(f_{\epsilon'}-f_{\epsilon})-\tilde{\lambda} \bar{g}_{\epsilon}^2)\\
%&=z(f_{\epsilon'}-f_{\epsilon})(g_{\epsilon}+\bar{g}_{\epsilon})-\tilde{\lambda} (g_{\epsilon}^2-\bar{g}_{\epsilon}^2)\\
&=z(f_{\epsilon'}-f_{\epsilon})(g_{\epsilon}+\bar{g}_{\epsilon})-\tilde{\lambda} (g_{\epsilon}^2-\bar{g}_{\epsilon}^2)\\
&=(g_{\epsilon}+\bar{g}_{\epsilon})(z(f_{\epsilon'}-f_{\epsilon})-\tilde{\lambda} (g_{\epsilon}-\bar{g}_{\epsilon}))
 \end{split}
\end{equation*}
and thus
$$z(f_{\epsilon'}-f_{\epsilon})-\tilde{\lambda} (g_{\epsilon}-\bar{g}_{\epsilon})=0,$$
which is the statement (iii).

Now we prove (iv). 
From (iii) we have 
\begin{equation*}
\tilde{\lambda} \bar{g}_{\epsilon}=\tilde{\lambda} g_{\epsilon}-z(f_{\epsilon'}-f_{\epsilon})
\end{equation*}
and thus  
$$1+\tilde{\lambda} \bar{g}_{\epsilon}g_{\epsilon}=1+(\tilde{\lambda} g_{\epsilon}-z(f_{\epsilon'}-f_{\epsilon}))g_{\epsilon}=1+\tilde{\lambda} g_{\epsilon}^2-zg_{\epsilon}(f_{\epsilon'}-f_{\epsilon})=f_{\epsilon} f_{\epsilon'}$$
by  (\ref{det-w}).
\end{proof}

%\subsection{Recursive formulas using the $\epsilon_i$-sequences }%=====================================
Now we consider the set of symmetric $ e_i$-sequences 
$$  {\bf e}(n)=( e_{n}, \cdots, e_2, e_1, e_1, e_2,\cdots, e_{n}),$$
for any sequence $ e_1,  e_2, e_3,\cdots$ with $ e_i \in\{-1, 1\} $, 
and denote $f_{ e(n)}$ and $ g_{ e(n)}$ by $f_n$ and $g_n$ so that $f_n$($=f_{\alpha-1}$ in Theorem \ref{kmot-rep2}) is the Riley polynomial of the kmot group $G(\epsilon)$, where
$$(\epsilon_1,\epsilon_2,\cdots,\epsilon_{\alpha-1})={\bf \epsilon}={\bf e}(n)=( e_{n}, \cdots, e_2, e_1, e_1, e_2,\cdots, e_{n}).$$
% Let $f_{0}(\lambda, \tilde{\lambda})=1,\,\,  g_{0}(M, \tilde{\lambda})=0.$
Then using the new $\epsilon$-notation  Proposition \ref{w-quandle} is restated as follows:
%\begin{equation}\label{w-n}
%W( {\bf e}(n))=\begin{pmatrix}
%f_{n}+zg_{n} & g_{n}\\
%  \lambda g_{n} & f_{n-1}-zg_{n}
% \end{pmatrix}
%\end{equation}  and  
%\begin{equation}\label{eq-2}
%1+\tilde{\lambda}g_{n}\bar{g}_{n}=f_nf_{n-1}
%\end{equation}
%\begin{equation}\label{eq-1}
%g_{n}= \bar{g}_{n-1}M^{-2 e_n}+ e_nf_{n-1}M^{- e_n},
%\end{equation}
%\begin{equation}\label{eq-3}
%f_{n}+zg_{n}= (f_{n-2}M^{2 e_n}+\lambda f_{n-1})+
%\bar{g}_{n-1}(2 e_n\lambda M^{ e_n}+z(-\lambda+M^{2 e_n}) )
%\end{equation} 
\begin{corollary}\label{cor-6.1}
 Let $f_{0}(\lambda, \tilde{\lambda})=1,\,\,  g_{0}(M, \tilde{\lambda})=0.$ Then 
\begin{enumerate}
   \item [\rm (i)] 
 $W( {\bf e}(n))=\begin{pmatrix}
f_{n}+zg_{n} & g_{n}\\
  \lambda g_{n} & f_{n-1}-zg_{n}
 \end{pmatrix}$  
  \item [\rm (ii)]  $1+\tilde{\lambda}g_{n}\bar{g}_{n}=f_nf_{n-1}$
   \item [\rm (iii)] $g_{n}= \bar{g}_{n-1}M^{-2 e_n}+ e_nf_{n-1}M^{- e_n}$
    \item [\rm (iv)]
   $ f_{n}+zg_{n}= (f_{n-2}M^{2 e_n}+\lambda f_{n-1})+
\bar{g}_{n-1}(2 e_n\lambda M^{ e_n}+z(-\lambda+M^{2 e_n}))$
\end{enumerate}
\end{corollary}
From this, we obtain the following Le's result (see Theorem 3.3.1 in \cite{Le}).
\begin{corollary}[\cite{Le}] The Riley polynomial $f_n$ of the kmot group $G( {\bf e}(n))$ satisfies the following:
\begin{enumerate}
    \item [\rm (i)] $f_n=(-1)^n+\sum_{k=1}^{n} (-1)^{n-k}tr\left(W( {\bf e}(k)\right)$
    \item [\rm (ii)]
    $\lambda f_n=tr\left(W( {\bf e}(n))\right)-tr\left(\begin{pmatrix}
M & 0\\
  \lambda& \frac{1}{M} 
 \end{pmatrix}
 W({\bf e}(n))
 \begin{pmatrix}
M & 1\\
  0& \frac{1}{M} 
 \end{pmatrix}^{-1}\right)$
\end{enumerate}
\end{corollary}
\begin{proof}
Since $tr\left(W( {\bf e}(k)\right)=f_k+f_{k-1}$ for each $k \in \mathbb N$ by Corollary \ref{cor-6.1}, (i) is proved as follows:
\begin{equation*}
    \begin{split}
        f_n&=(f_n+f_{n-1})-(f_{n-1}+f_{n-2})+(f_{n-2}+f_{n-3})-\cdots+(-1)^{n-1}(f_1+f_0)+(-1)^n\\
        &=tr\left(W( {\bf e}(n)\right)-tr\left(W( {\bf e}(n-1)\right)+\cdots+(-1)^{n-1}tr\left(W( {\bf e}(1)\right)+(-1)^n\\
        &=(-1)^n+\sum_{k=1}^{n} (-1)^{n-k}tr\left(W( {\bf e}(k)\right)
   \end{split}
   \end{equation*}
From  
\begin{equation*}
    \begin{split}
   \begin{pmatrix}
M & 0\\
  \lambda& \frac{1}{M} 
 \end{pmatrix}
 W({\bf e}(n))
 \begin{pmatrix}
M & 1\\
  0& \frac{1}{M} 
 \end{pmatrix}^{-1}&= \begin{pmatrix}
M & 0\\
  \lambda& \frac{1}{M} 
 \end{pmatrix}  \begin{pmatrix}
f_{n}+zg_{n} & g_{n}\\
  \lambda g_{n} & f_{n-1}-zg_{n}
 \end{pmatrix} 
  \begin{pmatrix}
\frac{1}{M} & -1\\
  0&  M
 \end{pmatrix}\\
 &= \begin{pmatrix}
M & 0\\
  \lambda& \frac{1}{M} 
 \end{pmatrix}
 \begin{pmatrix}
M^{-1}(f_{n}+zg_{n}) & -f_n+(M-z)g_{n}\\
  \lambda M^{-1} g_{n} & Mf_{n-1}-(\lambda+zM)g_{n}
 \end{pmatrix} \\
 &= \begin{pmatrix}
f_{n}+zg_{n} & **\\
  * & -\lambda f_n+f_{n-1}-zg_{n}
 \end{pmatrix} 
     \end{split}
   \end{equation*}
 we have $$tr\left(\begin{pmatrix}
M & 0\\
  \lambda& \frac{1}{M} 
 \end{pmatrix}
 W({\bf e}(n))
 \begin{pmatrix}
M & 1\\
  0& \frac{1}{M} 
 \end{pmatrix}^{-1}\right)=(f_n+f_{n-1})-\lambda f_n
 =tr\left(W( {\bf e}(n)\right)-\lambda f_n,
 $$
 which proves (ii).
\end{proof}
Now we have the following formulas for $f_n$ and $g_n$.
\begin{theorem}\label{f-g-formula}
Let ${\bf e}(n)=( e_{n}, \cdots, e_2, e_1, e_1, e_2,\cdots, e_{n})$. Then the polynomial $g_n=g_{{\bf e}(n)}$ and 
the Riley polynomial  $f_n=f_{{\bf e}(n)}$ of the kmot group $G({\bf e}(n))$ satisfy 
\begin{equation*}
\begin{split}
{\rm (i)}&
\left\{
\begin{array}{ll}
 f_{n}&= M^{2 e_n}f_{n-2}+(\tilde{\lambda}+1-M^{2 e_n}) f_{n-1}+ e_n\tilde{\lambda}(M+M^{-1})\bar{g}_{n-1} \\
g_{n}&=e_nM^{- e_n}f_{n-1}+M^{-2 e_n}\bar{g}_{n-1}
\end{array}			
\right.\\
{\text or} & \text{equivalently}\\
{\rm(ii)}& 
\left\{
\begin{array}{ll}
 f_{n}&=(\tilde{\lambda}+1-M^{2 e_n}) f_{n-1}+M^{2 e_n}f_{n-2}+ \tilde{\lambda}(M+M^{-1}) \sum_{i=1}^{n-1}e_ne_i M^{\mu_i}f_{i-1}\\
g_{n}&=\sum_{i=1}^{n}e_i M^{\mu_i}f_{i-1}\\
\end{array}			
\right.\\
%&=f_{n-2}M^{2 e_n}+(\lambda-1+M^{-2 e_n}) f_{n-1}+ e_n\tilde{\lambda}(M+M^{-1})( \displaystyle\sum_{i=1}^{n-1}e_i M^{(-1)^{n-i}e_i-2\sum_{j=i}^{n-1}(-1)^{n-j}e_j}f_{i-1})
% &\qquad\quad =e_nM^{- e_n}f_{n-1}+ e_{n-1}M^{-2 e_n+ e_{n-1}}f_{n-2}+\cdots+ e_{1}M^{-2 e_n+2 e_{n-1}+\cdots+(-1)^{n-1}2 e_{2}+(-1)^n e_{1}}\\
 \end{split}
 \end{equation*}
 where $\mu_i=(-1)^{n-i}e_i-2\sum_{j=i}^{n}(-1)^{n-j}e_j$.
\end{theorem}
\begin{proof}
(i)   The identity $g_{n}=e_nM^{- e_n}f_{n-1}+M^{-2e_n}\bar{g}_{n-1}$ is exactly (iii) of Corollary \ref{cor-6.1}, and by applying this identity starting from the formula
$$f_{n}= (f_{n-2}M^{2 e_n}+\lambda f_{n-1})+\bar{g}_{n-1}(2 e_n\lambda M^{ e_n}+z(-\lambda+M^{2 e_n})  )-zg_{n}\\$$
in Corollary \ref{cor-6.1}, we have the identity for $f_n$ as follows:
\begin{equation*}
\begin{split}
f_{n}&= f_{n-2}M^{2 e_n}+\lambda f_{n-1}+\bar{g}_{n-1}(2 e_n\lambda M^{ e_n}+z(-\lambda+M^{2 e_n})  )-zg_{n}\\
&= f_{n-2}M^{2 e_n}+\lambda f_{n-1}+\bar{g}_{n-1}(2 e_n\lambda M^{ e_n}+z(-\lambda+M^{2 e_n}))-z(\bar{g}_{n-1}M^{-2 e_n}+ e_nf_{n-1}M^{- e_n})\\
&= f_{n-2}M^{2 e_n}+(\lambda-z e_nM^{- e_n}) f_{n-1}+\bar{g}_{n-1}(2 e_n\lambda M^{ e_n}+z(-\lambda+M^{2 e_n}-M^{-2 e_n})  )\\
&=f_{n-2}M^{2 e_n}+(\lambda-1+M^{-2 e_n}) f_{n-1}+\bar{g}_{n-1}(2 e_n\lambda M^{ e_n}+z(-\lambda+M^{2 e_n}-M^{-2 e_n})  )\\
&= f_{n-2}M^{2 e_n}+(\tilde{\lambda}+1-M^{2 e_n})f_{n-1}+ e_n\bar{g}_{n-1}\tilde{\lambda}(M+M^{-1})\\
%&=  f_{n-2}M^{2 e_n}+(\lambda-1+M^{-2 e_n}) f_{n-1}+ e_n\tilde{\lambda}(M+M^{-1})( \displaystyle\sum_{i=1}^{n-1}e_i M^{(-1)^{n-i}e_i-2\sum_{j=i}^{n-1}(-1)^{n-j}e_j}f_{i-1})
 \end{split}
\end{equation*}
Here the last equality follows from
\begin{equation*}
\begin{split}
z(-\lambda+M^{2 e_n}-M^{-2 e_n})&=-z\lambda+ e_n(M^3+M^{-3}-(M+M^{-1}))\\
%&=-z\lambda+ e_n(M+M^{-1})(M^2+M^{-2}-2)\\
&=-z\lambda+ e_nz^2(M+M^{-1})\\
&=-z\lambda+ e_n(\tilde{\lambda}-\lambda)(M+M^{-1})\\
 \end{split}
\end{equation*}
and 
\begin{equation*}
\begin{split}
2 e_n\lambda M^{ e_n}+z(-\lambda+M^{2 e_n}-M^{-2 e_n})&=2 e_n\lambda M^{ e_n}-z\lambda+ e_n(\tilde{\lambda}-\lambda)(M+M^{-1})\\
&= e_n\tilde{\lambda}(M+M^{-1}).
 \end{split}
\end{equation*}

 (ii)   
By applying the identity for $g_n$ in (i) repeatedly, we obtain the formula for $g_n$ as follows:
\begin{equation*}
\begin{split}
g_{n}&=  e_nM^{- e_n}f_{n-1}+M^{-2 e_n}\bar{g}_{n-1}\\
&=  e_nM^{- e_n}f_{n-1}+M^{-2 e_n}( e_{n-1}M^{ e_{n-1}}f_{n-2}+M^{2 e_{n-1}}g_{n-2})\\
%&=  e_nM^{- e_n}f_{n-1}+ e_{n-1}M^{-2 e_n+ e_{n-1}}f_{n-2}+M^{-2 e_n+2 e_{n-1}}g_{n-2}\\
&=  e_nM^{- e_n}f_{n-1}+ e_{n-1}M^{-2 e_n+ e_{n-1}}f_{n-2}+M^{-2 e_n+2 e_{n-1}}( e_{n-2}M^{- e_{n-2}}f_{n-3}+M^{-2 e_{n-2}}\bar{g}_{n-3})\\
&=  \cdots\\
%&=  e_nM^{- e_n}f_{n-1}+ e_{n-1}M^{-2 e_n+ e_{n-1}}f_{n-2}+\cdots+ e_{1}M^{-2 e_n+2 e_{n-1}+\cdots+(-1)^{n-1}2 e_{2}+(-1)^n e_{1}}\\
&=\displaystyle\sum_{i=1}^{n}e_i M^{(-1)^{n-i}e_i-2\sum_{j=i}^{n}(-1)^{n-j}e_j}f_{i-1}
 \end{split}
\end{equation*}
Applying this identity to (i), we have the identity for $f_n$ immediately.
\end{proof}

\begin{example} By Corollary \ref{cor-6.1} or Theorem \ref{f-g-formula}, we have
\begin{equation*}
f_1=1+\tilde{\lambda},\,\,\,\,g_1= e_1M^{- e_1}f_{0}+M^{-2 e_1}\bar{g}_{0}= e_1M^{- e_1}\\
\end{equation*}
and
\begin{equation*}
\begin{split}
f_2&= (f_0M^{2 e_2}+(\lambda-1+M^{-2 e_2}) f_{1})+ e_2\bar{g}_{1}\tilde{\lambda}(M+M^{-1})\\
%&= (M^{2 e_2}+(\lambda-1+M^{-2 e_2}) (1+\tilde{\lambda}))+ e_2 e_1M^{ e_1}\tilde{\lambda}(M+M^{-1})\\
&=\left\{
\begin{array}{ll}
 (M^{2 e_2}+(\lambda-1+M^{-2 e_2}) (1+\tilde{\lambda}))+M^{ e_2}\tilde{\lambda}(M+M^{-1}) & \textrm{when} \quad  e_1= e_2 \\
 (M^{2 e_2}+(\lambda-1+M^{-2 e_2}) (1+\tilde{\lambda}))-M^{- e_2}\tilde{\lambda}(M+M^{-1})  & \textrm{when} \quad  e_1=- e_2 
\end{array}			
\right.\\
%&=\left\{
%\begin{array}{ll}
%\tilde{\lambda}(\lambda-1+M^{-2 e_2}+M^{1+ e_2}+M^{-1+ e_2}) + %(M^{2 e_2}+M^{-2 e_2}+\lambda-1)  & \textrm{when} \quad  e_1= e_2 \\
%\tilde{\lambda}(\lambda-1+M^{-2 e_2}-M^{1- e_2}-M^{-1- e_2}) + %(M^{2 e_2}+M^{-2 e_2}+\lambda-1)  & \textrm{when} \quad  e_1=- %e_2 
%\end{array}			
%\right.\\
%&=\left\{
%\begin{array}{ll}
%\tilde{\lambda}(\lambda+M^{2}+M^{-2})  + %(M^{2}+M^{-2}+\lambda-1)  & \textrm{when} \quad  e_1= e_2 \\
%\tilde{\lambda}(\lambda-2) + (M^{2}+M^{-2}+\lambda-1)  & %\textrm{when} \quad  e_1=- e_2 
%\end{array}			
%\right.\\
%&=\left\{
%\begin{array}{ll}
%\tilde{\lambda}(\lambda+z^2+2) + (z^2+\lambda+1)  & %\textrm{when} \quad  e_1= e_2 \\
%\tilde{\lambda}(\lambda-2) + (z^2+\lambda+1)  & \textrm{when} %\quad  e_1=- e_2 
%\end{array}			
%\right.\\
%&=\left\{
%\begin{array}{ll}
%\tilde{\lambda}(\tilde{\lambda}+2) + (\tilde{\lambda}+1)  & %\textrm{when} \quad  e_1= e_2 \\
%\tilde{\lambda}(\lambda-2) + (\tilde{\lambda}+1)  & %\textrm{when} \quad  e_1=- e_2 
%\end{array}			
%\right.\\
&=\left\{
\begin{array}{ll}
\tilde{\lambda}^2+3\tilde{\lambda}+1  & \textrm{when} \quad  e_1= e_2 \\
\tilde{\lambda}\lambda-\tilde{\lambda}+1 & \textrm{when} \quad  e_1=- e_2 
\end{array}			
\right.\\
g_2&= e_{2}M^{- e_2}f_1+ e_1M^{-2 e_2+ e_1}= e_{2}M^{- e_2}(1+\tilde{\lambda})+ e_1M^{-2 e_2+ e_1}\\
&=\left\{
\begin{array}{ll}
 e M^{- e}(1+\tilde{\lambda})+ e M^{- e}= e M^{- e}(2+\tilde{\lambda})  & \textrm{when} \quad  e_1= e_2= e \\
 e M^{- e}(1+\tilde{\lambda})- e M^{-3 e}=z e M^{-2 e}+M^{- e}\tilde{\lambda}  & \textrm{when} \quad  e_1=- e_2= e 
\end{array}			
\right.\\
\end{split}
\end{equation*}
\end{example}

%\subsection{Recursive formulas using the sign change sequences }%=====================================
By (\ref{f-g}), we know that the Riley polynomials depend only on the sign changes, which implies that the Riley polynomials for the kmot group $G(\epsilon)$ and its ``mirror" $G(-\epsilon)$ are the same each other. The following recursive formula also shows the property.
\begin{corollary} \label{delta-formula}
Let $\delta_n=e_{n-1}e_n$. Then 
 the following recursive formula holds for all $n>2$.
\begin{equation*}
\begin{split}
f_n&=\left\{
\begin{array}{ll}
 (\tilde{\lambda}+ 1)(f_{n-1}+f_{n-2})-f_{n-3}  & \textrm{when} \quad  \delta_n=1 \\
  (\lambda-1)(f_{n-1}-f_{n-2})+f_{n-3}   & \textrm{when} \quad   \delta_n=-1
\end{array}			
\right.\\
&=(\tilde{\lambda}+\frac{\delta_n-1}{2}z^2+\delta_n)(f_{n-1}+\delta_nf_{n-2})-\delta_nf_{n-3}
\end{split}
\end{equation*}
%where $A=\tilde{\lambda}+\frac{e_{n-1}e_n-1}{2}s+e_{n-1}e_n$.
\end{corollary}
\begin{proof}
By Corollary \ref{cor-6.1}, (iii) and Theorem \ref{f-g-formula}, (ii),  
\begin{equation*}
\begin{split}
 f_{n}&= f_{n-2}M^{2 e_n}+(\lambda-1+M^{-2 e_n}) f_{n-1}+ e_n\bar{g}_{n-1}\tilde{\lambda}(M+M^{-1})\\
&= f_{n-2}M^{2 e_n}+(\lambda-1+M^{-2 e_n}) f_{n-1}+ e_n\tilde{\lambda}(M+M^{-1})[g_{n-2}M^{2 e_{n-1}}+ e_{n-1}f_{n-2}M^{ e_{n-1}}]\\
\end{split}
\end{equation*}
and 
$$\tilde{\lambda}(M+M^{-1})g_{n-2}=e_{n-1}[f_{n-1}-f_{n-3}M^{-2e_{n-1}}-(\lambda-1+M^{2 e_{n-1}}) f_{n-2}].$$
Hence we have
\begin{equation*}
\begin{split}
 f_{n}&=(\lambda-1+M^{-2 e_n}) f_{n-1}+ [M^{2 e_n}+e_{n-1}e_n\tilde{\lambda}(M+M^{-1})M^{ e_{n-1}}]f_{n-2}\\
&\quad+e_{n-1}e_nM^{2 e_{n-1}}[f_{n-1}-f_{n-3}M^{-2e_{n-1}}-(\lambda-1+M^{2 e_{n-1}}) f_{n-2}]\\
&=A f_{n-1}+Bf_{n-2}- \delta_nf_{n-3}\\
\end{split}
\end{equation*}
where 
\begin{equation*}
\begin{split}
A&=\lambda-1+M^{-2 e_n}+\delta_nM^{2 e_{n-1}}\\
B&=M^{2 e_n}+\delta_n\tilde{\lambda}(M+M^{-1})M^{ e_{n-1}}-\delta_nM^{2 e_{n-1}}(\lambda-1+M^{2 e_{n-1}}).    
\end{split}
\end{equation*}

So it suffices to show that 
\begin{equation*}
\begin{split}
A&=\left\{
\begin{array}{ll}
 \tilde{\lambda}+ 1  & \textrm{when} \quad  \delta_n=1 \\
  \lambda-1 & \textrm{when} \quad   \delta_n=-1
\end{array}			
\right.\\
B&=\delta_nA=\left\{
\begin{array}{ll}
 \tilde{\lambda}+ 1  & \textrm{when} \quad  \delta_n=1 \\
  -\lambda+1 & \textrm{when} \quad   \delta_n=-1
\end{array}			
\right.\\
\end{split}
\end{equation*}
since $\tilde{\lambda}+\frac{\delta_n-1}{2}z^2+\delta_n=
\left\{
\begin{array}{ll}
 \tilde{\lambda}+ 1  & \textrm{when} \quad  \delta_n=1 \\
  -\lambda+1 & \textrm{when} \quad   \delta_n=-1
\end{array}			
\right.$.

The first identity is proved as follows.
\begin{equation*}
\begin{split}
A&=\lambda-1+M^{-2 e_n}+\delta_nM^{2 e_{n-1}}\\
&=\left\{
\begin{array}{ll}
 \lambda-1+M^{-2 e_{n-1}}+M^{2 e_{n-1}}  & \textrm{when} \quad  \delta_n=1 \\
 \lambda-1+M^{2 e_{n-1}}-M^{2 e_{n-1}}  & \textrm{when} \quad   \delta_n=-1
\end{array}			
\right.\\
%&=\left\{
%\begin{array}{ll}
% \lambda-1+M^{2}+M^{-2}  & \textrm{when} \quad  \delta_n=1 \\
% \lambda-1  & \textrm{when} \quad   \delta_n=-1
%\end{array}			
%\right.\\
%&=\left\{
%\begin{array}{ll}
% \lambda+z^2+1  & \textrm{when} \quad  \delta_n=1 \\
% \lambda-1  & \textrm{when} \quad   \delta_n=-1
%\end{array}			
%\right.\\
&=\left\{
\begin{array}{ll}
 \tilde{\lambda}+ 1  & \textrm{when} \quad  \delta_n=1 \\
  \lambda-1 & \textrm{when} \quad   \delta_n=-1
\end{array}			
\right.\\
\end{split}
\end{equation*}
The second identity is proved as follows. If $\delta_n=1$, that is, $e_{n-1}=e_n $, then 
\begin{equation*}
\begin{split}
B&=M^{2 e_n}+\delta_n\tilde{\lambda}(M+M^{-1})M^{ e_{n-1}}-\delta_nM^{2 e_{n-1}}(\lambda-1+M^{2 e_{n-1}})\\
&=2M^{2 e_n}+\tilde{\lambda}(M+M^{-1})M^{ e_{n}}-M^{2 e_{n}}(\lambda+M^{2 e_{n}})  \\
%&= \tilde{\lambda}+(\tilde{\lambda}-\lambda-2) M^{2 e_{n}}-M^{4 e_{n}}   \\
%&= \tilde{\lambda}+(z^2+2-M^{2 e_{n}} )M^{2 e_{n}}  \\
%&= \tilde{\lambda}+(M^{-2 e_{n}} )M^{2 e_{n}}  \\
&= \tilde{\lambda}+ 1 \\
&=\delta_nA
\end{split}
\end{equation*}
and if $\delta_n=-1$, that is, $e_{n-1}=-e_n$, then 
\begin{equation*}
\begin{split}
B&=M^{2 e_n}+\delta_n\tilde{\lambda}(M+M^{-1})M^{ e_{n-1}}-\delta_nM^{2 e_{n-1}}(\lambda-1+M^{2 e_{n-1}})\\
&= M^{2 e_n}-\tilde{\lambda}(M+M^{-1})M^{-e_{n}}+M^{-2e_{n}}(\lambda-1+M^{-2e_{n}}) \\
%&= M^{2 e_n}-\tilde{\lambda}+M^{-2e_{n}}(\lambda-\tilde{\lambda}-1+M^{-2e_{n}}) \\
%&=  M^{2 e_n}-\tilde{\lambda}+M^{-2e_{n}}(-z^2-1+M^{-2e_{n}})  \\
%&=  M^{2 e_n}-\tilde{\lambda}+M^{-2e_{n}}(1-M^{2e_{n}})  \\
%&=  M^{2 e_n}+M^{-2e_{n}}-1-\tilde{\lambda}  \\
%&=  z^2+1-\tilde{\lambda}  \\
&= -\lambda+ 1 \\
&=\delta_nA
\end{split}
\end{equation*}
\end{proof}

%\begin{remark}
Note that the recursive formula in Corollary \ref{delta-formula} is equivalent to the ones in Lemma 1.2 of Burde \cite{Burde} and Proposition 5.2 of Hilden-Lozano-Montesinos \cite{HLM}.
  Now we prove the following recursive formula for Alexander polynomials using Corollary \ref{delta-formula}.
%\end{remark}

\begin{corollary}\label{alexander2}
Let $\bigtriangleup_n(t)$ be the Alexander polynomial of the kmot group $G({\bf e}(n))$ and  $\alpha_n(t)=\frac{\delta_n+1}{2}(t+t^{-1})$. Then
%$$ \bigtriangleup_n(t)=(\alpha_n(t)-1)\bigtriangleup_{n-1}(t)+(\alpha_n(t)-\delta_n)\bigtriangleup_{n-2}(t)-\delta_n\bigtriangleup_{n-3}(t)$$
$$ \bigtriangleup_n(t)+\bigtriangleup_{n-1}(t)=\alpha_n(t)(\bigtriangleup_{n-1}(t)+\bigtriangleup_{n-2}(t))-\delta_n(\bigtriangleup_{n-2}(t)+\bigtriangleup_{n-3}(t))$$
\end{corollary}
\begin{proof}
For the case of $\lambda=0$, $\tilde{\lambda}$ is equal to $z^2$ and 
\begin{equation*}
\begin{split}
f_n&=(\tilde{\lambda}+\frac{\delta_n-1}{2}z^2+\delta_n)(f_{n-1}+\delta_nf_{n-2})-\delta_nf_{n-3}\\
&=(z^2+\frac{\delta_n-1}{2}z^2+\delta_n)(f_{n-1}+\delta_nf_{n-2})-\delta_nf_{n-3}\\
&=(\frac{\delta_n+1}{2}z^2+\delta_n)(f_{n-1}+\delta_nf_{n-2})-\delta_nf_{n-3}\\
&=(\frac{\delta_n+1}{2}z^2+\delta_n)f_{n-1}+(\frac{\delta_n+1}{2}z^2+1)f_{n-2}-\delta_nf_{n-3}\\
&=(\frac{\delta_n+1}{2}(M^2+M^{-2})-1)f_{n-1}+(\frac{\delta_n+1}{2}(M^2+M^{-2})-\delta_n)f_{n-2}-\delta_nf_{n-3}\\
&=(\alpha_n(M^2)-1)f_{n-1}+(\alpha_n(M^2)-\delta_n)f_{n-2}-\delta_nf_{n-3}.
\end{split}
\end{equation*}
So we have 
$$ \bigtriangleup_n(M^2)=(\alpha_n(M^2)-1)\bigtriangleup_{n-1}(M^2)+(\alpha_n(M^2)-\delta_n)\bigtriangleup_{n-2}(M^2)-\delta_n\bigtriangleup_{n-3}(M^2,)$$
which implies 
$$ \bigtriangleup_n(t)=(\alpha_n(t)-1)\bigtriangleup_{n-1}(t)+(\alpha_n(t)-\delta_n)\bigtriangleup_{n-2}(t)-\delta_n\bigtriangleup_{n-3}(t)$$
or equivalently
$$ \bigtriangleup_n(t)+\bigtriangleup_{n-1}(t)=\alpha_n(t)(\bigtriangleup_{n-1}(t)+\bigtriangleup_{n-2}(t))-\delta_n(\bigtriangleup_{n-2}(t)+\bigtriangleup_{n-3}(t)).$$
\end{proof}
From the recursive formula in Corollary \ref{alexander2}, we can prove that $\bigtriangleup_n(t)+\bigtriangleup_{n-1}(t)$ is actually a well-known Chebyshev polynomial. Here, $Chebyshev$ $polynomials$ are polynomials defined by a three-term recursion
$$g_{n+1}(x)=xg_n(x)-g_{n-1}(x).$$ 
 We denote the Chebyshev polynomials $g_n(x)$ with the intial condition 
$g_0(x)=a, g_1(x)=b$ by $Ch_n^x(a,b)$.
The following Chebyshev polynomials are frequently used: 
\begin{equation*}
\begin{split}
s_n(x)&=Ch_n^x(0,1)\\
v_n(x)&=Ch_n^x(1,x-1)=s_{n+1}(x)-s_n(x)\\
t_n(x)&=Ch_n^x(2,x)=s_{n+1}(x)-s_{n-1}(x)
\end{split}
\end{equation*}
Note that $s_{n-1}(x)=-Ch_n^x(1,0)$ and 
$$Ch_n^x(a,b)=aCh_n^x(1,0)+bCh_n^x(0,1)=-as_{n-1}(x)+bs_n(x).$$
Note that $T_n(x)=\frac{1}{2}t_{n}(2x)$, $U_n(x)=s_{n+1}(2x)$, and $V_n(x)=v_n(2x)$ are exactly the Chebyshev polynomials of the first, second, and third  kinds, respectively. 

We have the following identity.
\begin{proposition}\label{alexaner3}
Let $\bigtriangleup_n(t)$ be the Alexander polynomial of the kmot group $G({\bf e}(n))$ and   
$\beta_n=\sum_{i=1}^n \delta_1\cdots\delta_i$.(We let $\delta_1=1$.)  Then
$$ \bigtriangleup_{n+1}(t)+\bigtriangleup_{n}(t)=t_{\beta_{n+1}}(t+t^{-1})=t^{\beta_{n+1}}+t^{-\beta_{n+1}}.$$
%where $t_m(x)=s_{m+1}(x)-s_{m-1}(x)=Ch_m^x(2,x)$. 
\end{proposition}
\begin{proof}
It suffices to show the first equality 
since $ t_m(x+x^{-1})=x^m+x^{-m}$ is an well-known fact.
Let $Q_n(t)=\bigtriangleup_n(t)+\bigtriangleup_{n-1}(t)$. Then 
$$Q_{n+1}(t)=\left\{
\begin{array}{ll}
 (t+t^{-1})Q_n(t)-Q_{n-1}(t)  & \textrm{when} \quad  \delta_{n+1}=1 \\
  Q_{n-1}(t) & \textrm{when} \quad   \delta_{n+1}=-1
\end{array}			
\right.$$ by Corollary \ref{alexander2}, and 
\begin{equation*}
\begin{split}
Q_1(t)&=\bigtriangleup_1(t)+\bigtriangleup_{0}(t)=(1+s)+1=2+s=t+t^{-1}=t_1(t+t^{-1})=t_{\beta_1}(t+t^{-1})\\
Q_2(t)&=\bigtriangleup_2(t)+\bigtriangleup_{1}(t)\\
&=\left\{
\begin{array}{ll}
 (s^3+3s+1)+(1+s) & \textrm{when} \quad  \delta_2=1 \\
  (1-s)+(1+s) & \textrm{when} \quad   \delta_2=-1
\end{array}			
\right.\\
&=\left\{
\begin{array}{ll}
 (s+2)^2-2  & \textrm{when} \quad  \delta_2=1 \\
  2 & \textrm{when} \quad   \delta_2=-1
\end{array}			
\right.\\
&=\left\{
\begin{array}{ll}
 (t+t^{-1})^2-2  & \textrm{when} \quad  \delta_2=1 \\
  2 & \textrm{when} \quad   \delta_2=-1
\end{array}			
\right.\\
&=\left\{
\begin{array}{ll}t_2(t+t^{-1})  & \textrm{when} \quad  \delta_2=1 \\
  t_0(t+t^{-1}) & \textrm{when} \quad   \delta_2=-1
\end{array}			
\right.\\
&=t_{\beta_2}(t+t^{-1})
\end{split}
\end{equation*}
since $\beta_1=\delta_1=1$ and 
$\beta_1=\left\{
\begin{array}{ll}
2  & \textrm{when} \quad  \delta_2=1 \\
  0 & \textrm{when} \quad   \delta_2=-1
\end{array}.
\right.$
Suppose that the statement holds for all $n\leq k$ to proceed by induction on $n$.
Then we have 
\begin{equation*}
\begin{split}
Q_{k+1}(t)&=\left\{
\begin{array}{ll}
 (t+t^{-1})Q_k(t)-Q_{k-1}(t)  & \textrm{when} \quad  \delta_{k+2}=1 \\
  Q_{k-1}(t) & \textrm{when} \quad   \delta_{k+2}=-1
\end{array}			
\right.\\
&=\left\{
\begin{array}{ll}
 (t+t^{-1})t_{\beta_{k+1}}(t+t^{-1})-t_{\beta_{k}}(t+t^{-1})& \textrm{when} \quad  \delta_{k+2}=1 \\
  t_{\beta_{k}}(t+t^{-1})& \textrm{when} \quad   \delta_{k+2}=-1
\end{array}.			
\right.
\end{split}
\end{equation*}
On the other hand, 
if $ \delta_{k+2}=-1$ then $$\beta_{k+2}=\beta_k+\delta_1\cdots\delta_{k+1}+\delta_1\cdots\delta_{k+1}\delta_{k+2}=\beta_k,$$ which implies $ t_{\beta_{k+2}}(t+t^{-1})= t_{\beta_{k}}(t+t^{-1})$, and if $\delta_{k+2}=1$ then 
$$\beta_{k+2}=\beta_k+\delta_1\cdots\delta_{k+1}+\delta_1\cdots\delta_{k+1}\delta_{k+2}=\beta_k+2\delta_1\cdots\delta_{k+1},$$ which implies 
 $ t_{\beta_{k+2}}(t+t^{-1})= (t+t^{-1}) t_{\beta_{k+1}}(t+t^{-1})-t_{\beta_{k}}(t+t^{-1})$
 since
 $\beta_{k}, \beta_{k+1}, \beta_{k+2}$ are consecutive numbers in this case.
 Therefore 
$$Q_{k+1}(t)=t_{\beta_{k+2}}(t+t^{-1})$$ 
in either case. This completes the proof.
\end{proof}

Note that $\beta_i=e_1(e_1+e_2+\cdots+e_i)$ for each $i=1,2,\cdots,n$ and thus $$t^{\beta_i}+t^{-\beta_i}=t^{e_1(e_1+e_2+\cdots+e_i)}+t^{-e_1(e_1+e_2+\cdots+e_i)}=t^{e_1+e_2+\cdots+e_i}+t^{-(e_1+e_2+\cdots+e_i)}.$$ So we have the following as a corollary of Proposition \ref{alexaner3}:
\begin{corollary}
$$\bigtriangleup_n(t)= (-1)^n+\sum_{i=1}^{n} (-1)^{n-i}(t^{\beta_i}+t^{-\beta_i})$$
or equivalently
$$(-1)^n\bigtriangleup_n(t)=1+\sum_{i=1}^{n} (-1)^{i}(t^{e_1+e_2+\cdots+e_i}+t^{-(e_1+e_2+\cdots+e_i)})$$
\end{corollary}
This result shows that the Alexander polynomial of the kmot group $G({\bf e}(n))$ is also the sum of Chebyshev polynomials: $\bigtriangleup_n(t)= (-1)^n+\sum_{i=1}^{n} (-1)^{n-i}t_{\beta_i}(t+t^{-1})$.
Furthermore we can obtain Minkus' formula again from the above formula as follows: if we let
$$ {\bf \epsilon}={\bf e}(n)=( e_{n}, \cdots, e_2, e_1, e_1, e_2,\cdots, e_{n}),$$ then 
$$\epsilon_1=e_n=\epsilon_{2n},\,\, \epsilon_2=e_{n-1}=\epsilon_{2n-1},\cdots, \epsilon_n=e_1=\epsilon_{n+1},$$ and thus
\begin{equation*}
\begin{split}
    t^{\frac{\sigma}{2}}\bigtriangleup_n(t)&= (-1)^nt^{\frac{\sigma}{2}}+t^{\frac{\sigma}{2}}\sum_{i=1}^{n} (-1)^{n-i}(t^{\beta_i}+t^{-\beta_i})\\
    &=(-1)^nt^{e_1+e_2+\cdots+e_n}+t^{e_1+e_2+\cdots+e_n}\sum_{i=1}^{n} (-1)^{n-i}(t^{e_1+e_2+\cdots+e_i}+t^{-(e_1+e_2+\cdots+e_i)})\\
    &=(-1)^nt^{e_1+e_2+\cdots+e_n}+\sum_{i=1}^{n} (-1)^{n-i}(t^{e_n+e_{n-1}+\cdots+e_1+e_1+e_2+\cdots+e_i}+t^{e_{n}+e_{n-1}+\cdots+e_{i+1})}\\
    &=(-1)^nt^{\epsilon_1+\epsilon_2+\cdots+\epsilon_{n}}+\sum_{i=1}^{n} (-1)^{n-i}(t^{\epsilon_1+\epsilon_2+\cdots+\epsilon_{n+i}}+t^{\epsilon_1+\epsilon_2+\cdots+\epsilon_{n-i}})\\
    &=1+\sum_{i=1}^{n}(-1)^{i}t^{\epsilon_1+\epsilon_2+\cdots+\epsilon_{i}}.
\end{split}
\end{equation*}

\section{2-bridge knots}%====================

For a knot $K$, the fundamental group of the complement  is called the $knot$ $group$  and denoted by $G(K)$. 
The knot group of  a 2-bridge knot $K=S(\alpha,\beta)$ always has a presentation of the form 
\begin{equation*}
G(K)=\pi_1(\Sphere^3 \setminus K)=\langle x,y \,|\,wx=yw \rangle,
\end{equation*}
where $w$ is of the form
\begin{equation*}
w=x^{\epsilon_1}y^{\epsilon_2}x^{\epsilon_3}y^{\epsilon_4}\cdots x^{\epsilon_{\alpha-2}}y^{\epsilon_{\alpha-1}}
\end{equation*}
with each $\epsilon_i=(-1)^{[i\frac{\beta}{\alpha}]}$. Note that  $\epsilon_i=\epsilon_{\alpha-i}$ for $i=1, \cdots, \alpha-1$ and thus $G(K)$ is a kmot group.  But there is a kmot group which is not a 2-bridge knot group, for example, 
a word $w=xy^{-1}x^{-1}y^{-1}x^{-1}y$ does not define any 2-bridge knot group, since for $\alpha=7$, $\{\pm(1,1,1,1,1,1), \pm(1,1,-1,-1,1,1), \pm(1,-1,1,1,-1,1)\}$ is the set of $\epsilon_i$-sequences from 2-bridge knot groups.
Note that some meridian $\mu$ and the preferred longitude $\gamma$ are represented by the words
$$\mu=x\quad\text{and}\quad \gamma=w^*wx^{-2\sigma}$$
where $\sigma=\displaystyle \sum_{i=1}^{\alpha-1} \epsilon_i$ and 
$
w^*=y^{\epsilon_1}x^{\epsilon_2}y^{\epsilon_3}x^{\epsilon_4}\cdots y^{\epsilon_{\alpha-2}}x^{\epsilon_{\alpha-1}}.
$

%\subsection{$SL(2,\bc)$-representations of 2-bridge knot groups} %------------------------------------------
Suppose that $\rho : G(K)  \rightarrow SL(2,\bc)$ is a non-abelian $ SL(2,\bc)$-representation. Then after conjugating if necessary, we may assume 
\begin{equation*}
\rho(x)=m=\begin{pmatrix}
M & 1\\
  0 & \frac{1}{M} 
 \end{pmatrix} \quad \text{and}\quad
\rho(y)=\begin{pmatrix}
M & 0\\
  \lambda& \frac{1}{M} 
 \end{pmatrix}.
\end{equation*}

Since $G(K)$ is a kmot group, 
a pair $M, \lambda$ determines a non-abelian representation if and only if  $W_{11}-zW_{12}=0$ by Proposition \ref{Riley poly}, where $W_{ij}$  is the ($i,j$)-element of 
\begin{equation*}
\begin{split}
W&=\rho(w)=\rho(x)^{\epsilon_1}\rho(y)^{\epsilon_2}\rho(x)^{\epsilon_3}\rho(y)^{\epsilon_4}\cdots \rho(x)^{\epsilon_{\alpha-2}}\rho(y)^{\epsilon_{\alpha-1}}.
\end{split}
\end{equation*}
Hence  $W_{11}-zW_{12}$ is called the $Riley$ $polynomial$ of $K$.

%Let  $\rho(w^*)=W^*$ and  $\tilde{W}(M,\lambda)=W(M^{-1},\lambda)$.

It is well-known that $\rho$-image of the longitude is 
\begin{equation}\label{longitude}
\rho(\gamma)=W^{*}Wm^{-2\sigma}=\begin{pmatrix}
L & *\\
 0& \frac{1}{L} 
 \end{pmatrix}, \quad  LM^{2\sigma}=-\frac{W_{12}}{\bar{W}_{12}}=\frac{W_{11}}{\bar{W}_{11}},
\end{equation}
which  can also be shown as follows: Since 
\begin{equation*}
\begin{split}
W^{*}Wm^{-2\sigma}&=\begin{pmatrix}
\bar{W}_{22} & \bar{W}_{12} \\
\lambda \bar{W}_{12}& \bar{W}_{11}
 \end{pmatrix} \begin{pmatrix}
W_{11} & W_{12} \\
\lambda W_{12}& W_{22}
 \end{pmatrix}
\begin{pmatrix}
M & 1\\
  0 & \frac{1}{M} 
 \end{pmatrix}^{-2\sigma}\\
&=\begin{pmatrix}
(W_{11}\bar{W}_{22}+\lambda W_{12}\bar{W}_{12})M^ {-2\sigma}& * \\
\lambda(W_{11}\bar{W}_{12}+\bar{W}_{11}W_{12})M^ {-2\sigma}& *
 \end{pmatrix}
\end{split}
\end{equation*}
by Lemma \ref{w-star},
\begin{equation}\label{w-equation}
W_{11}\bar{W}_{12}+\bar{W}_{11}W_{12}=(zW_{12})\bar{W}_{12}+(-z\bar{W}_{12})W_{12}=0
\end{equation}
and 
\begin{equation*}
\begin{split}
\bar{W}_{12}L&=\bar{W}_{12}(W_{11}\bar{W}_{22}+\lambda W_{12}\bar{W}_{12})M^ {-2\sigma}\\
&=(W_{11}\bar{W}_{12}\bar{W}_{22}+\lambda W_{12}\bar{W}_{12}^2)M^ {-2\sigma}\\
&=(-W_{12}\bar{W}_{11}\bar{W}_{22}+\lambda W_{12}\bar{W}_{12}^2)M^ {-2\sigma}  \quad\quad\text{using (\ref{w-equation})}\\   
&=(-W_{12}(1+\lambda\bar{W}_{12}^2\bar{W}_{22})+\lambda W_{12}\bar{W}_{12}^2)M^ {-2\sigma} \quad\quad\text{using } \det W^*=1\\
&=-W_{12}M^ {-2\sigma}.
\end{split}
\end{equation*}
(See \cite{HS} or \cite{Riley2}.)
So we have the following from (\ref{f-g}), (\ref{longitude}), Theorem \ref{kmot-rep2}, and Proposition \ref{w-quandle}, (iv).
\begin{corollary}\label{f-longitde}
Let $K=S(\alpha,\beta)$, $G=G(K)$, and $\epsilon_i=(-1)^{[i\frac{\beta}{\alpha}]}$ for $i=1, \cdots, \alpha-1$. Then
\begin{enumerate}
\item [\rm (i)]
$f_{\epsilon}(M,\tilde{\lambda})$ is the Riley polynomial of $K$,
\item [\rm (ii)]
$\displaystyle LM^{2\sigma}= -\frac{g_{\epsilon}(M,\tilde{\lambda})}{g_{\epsilon}(M^{-1},\tilde{\lambda})}=\frac{\tilde{\lambda} g_{\epsilon}(M,\tilde{\lambda})^2}{1-f_{\epsilon}(M,\tilde{\lambda})f_{\epsilon'}(M,\tilde{\lambda})}$,
\item [\rm (iii)]
$\displaystyle LM^{2\sigma}\equiv \tilde{\lambda}g_{\epsilon}(M,\tilde{\lambda})^2=\tilde{\lambda}g_{-\epsilon}(M^{-1},\tilde{\lambda})^2\quad ({\rm mod}\,\, f_{\epsilon}(M,\tilde{\lambda}))$.
\end{enumerate}
\end{corollary}
 By Corollary \ref{f-longitde}, the A-polynomial of a knot  $K=S(\alpha,\beta)$ is determined by two polynomials, $f_{\epsilon}(M,\tilde{\lambda})$ and $ LM^{2\sigma}g_{\epsilon}(M^{-1},\tilde{\lambda})+ g_{\epsilon}(M,\tilde{\lambda})$ (or  $LM^{2\sigma}-\tilde{\lambda}g_{\epsilon}(M,\tilde{\lambda})^2$), eliminating the variable  $\tilde{\lambda}$. This can be done by computing the Resultant of these two polynomials. Note that the resultant gives the A-polynomial of $K$ up to multiplicity, since $f_{\epsilon}(M,\tilde{\lambda})$ is monic in $\tilde{\lambda}$-variable. (See \cite{CCGLS, HS3,Huynh-Le, Petersen} for the definition and some basic properties of A-polynomials.) 
 
A-polynomials are known for most knots up to nine crossings, many knots with 10 crossings, several infinite families of 2-bridge knots, some pretzel knots, and all torus knots (refer to  Knot Table \cite{Rolfsen}, Knotinfo \cite{CL}, Culler list \cite{Culler}, \cite{HS} and \cite{Petersen}.) However, computing A-polynomials has been a challenging task.
 
Fortunately, our recursive formulas for $f_{\epsilon}$ and $g_{\epsilon}$ presented in Theorem \ref{f-g-formula} have proven to be quite effective in computing the A-polynomials of 2-bridge knots. Using Mathematica, we were able to obtain the A-polynomials for all 2-bridge knots with up to 12 crossings. Each computation took no more than a few minutes.

To provide some perspective on the computational efficiency, the average time required to compute the A-polynomial for a 10-crossing knot was approximately 12 seconds. For 11-crossing knots, it took around 102 seconds, and for 12-crossing knots, the computation time was approximately 1266 seconds.
We would like to invite interested readers to visit our website at https://2bridge.diagram.site to access the complete list of A-polynomials for 2-bridge knots with up to 12 crossings.
% These contain many unknown ones, that is, the  A-polynomials of  $10_{23}, 10_{25}, 10_{26},\cdots, 10_{45},11a_{13},11a_{59},11a_{65},\cdots$ are included. 
 %One can see  the full list of A-polynomials of 2-bridge knots of up to 12 crossings at https://2bridge.diagram.site. 
 
%Note that  $f_{\alpha-1}(M,\tilde{\lambda})=0$ implies that 
%\begin{equation*}
%\tilde{W}=
%\begin{pmatrix}
%zW_{12} & 0\\
%\lambda W_{12}& \lambda W_{12}-zW_{22}
% \end{pmatrix}\sim \begin{pmatrix}
%z & 0\\
%\lambda& -1
% \end{pmatrix}
%\sim\begin{pmatrix}
%zg_{\alpha-1}(M,\tilde{\lambda})  & 0\\
%\lambda g_{\alpha-1}(M,\tilde{\lambda})& \tilde{\lambda}g_{\alpha-1}(M^{-1},\tilde{\lambda}) 
% \end{pmatrix}
%=\tilde{c}_{\alpha-1}
%\end{equation*}
\section*{Acknowledgements }
The first author acknowledges support from the National Research Foundation of Korea (NRF) grant funded by the Korean government (MSIT) (NRF-2021R1F1A1049444) and the Open KIAS Center at the Korea Institute for Advanced Study.
The second author acknowledges support from the National Research Foundation of Korea (NRF) grant funded by the Korea government (MSIT) (NRF-2018R1A2B6005691).
The authors would like to express their gratitude to Philip Choi for his valuable assistance with the Mathematica computations.

\end{document}